\documentclass[12pt]{elsarticle}
\usepackage[margin = 2.5cm] {geometry} 

\usepackage{graphicx}
\usepackage{hyperref}		
\usepackage{amsmath}		
\usepackage{float}			
\usepackage[T1]{fontenc}    

\hypersetup {colorlinks, citecolor = blue, filecolor = black, linkcolor = blue, urlcolor = blue}

\begin{document}

\pagestyle{myheadings}
\thispagestyle{plain}
\markboth{KAROL MIKULA ET AL.}{NATURAL NUMERICAL NETWORKS FOR NATURA 2000 CLASSIFICATION}

\begin{frontmatter}

\title{Natural Numerical Networks for Natura 2000 habitats \\ classification by satellite images}

\author[1,2]{Karol Mikula\corref{cor1}}
\cortext[cor1]{{Corresponding author: Department of Mathematics, Slovak University of Technology in Bratislava, Radlinsk\'{e}ho 11, 810 05 Bratislava, Slovakia and Algoritmy:SK, s.r.o., \v{S}ulekova 6, 811 06 Bratislava, Slovakia.\\ E-mail address: {\tt mikula@math.sk}.}}
\author[1,2]{Michal Koll\'{a}r}
\author[1,2]{Aneta~A. O\v{z}vat}
\author[1,2]{Martin Ambroz}
\author[3]{Lucia \v{C}ahojov\'{a}}
\author[3]{Ivan Jarol\'{i}mek}
\author[3]{Jozef \v{S}ib\'{i}k}
\author[3]{M\'{a}ria \v{S}ib\'{i}kov\'{a}}

\address[1]{Department of Mathematics, Slovak University of Technology in Bratislava, Radlinsk\'{e}ho 11, 810 05 Bratislava, Slovakia.}
\address[2]{Algoritmy:SK, s.r.o., \v{S}ulekova 6, 811 06 Bratislava, Slovakia.}
\address[3]{Plant Science and Biodiversity Center, Slovak Academy of Sciences, Institute of Botany, D\'{u}bravsk\'{a} cesta 9, 845 23 Bratislava, Slovakia.}

\begin{abstract}
Natural numerical networks are introduced as a new classification algorithm based on the numerical solution of nonlinear partial differential equations of forward-backward diffusion type on complete graphs. The proposed natural numerical network is applied to open important environmental and nature conservation task, the automated identification of protected habitats by using satellite images. In the natural numerical network, the forward diffusion causes the movement of points in a feature space toward each other. The opposite effect, keeping the points away from each other, is caused by backward diffusion. This yields the desired classification. The natural numerical network contains a few parameters that are optimized in the learning phase of the method. After learning parameters and optimizing the topology of the network graph, classification necessary for habitat identification is performed. A relevancy map for each habitat is introduced as a tool for validating the classification and finding new Natura 2000 habitat appearances.

\end{abstract}

\begin{keyword}
data classification, partial differential equations on graphs, forward-backward diffusion, numerical methods, Natura 2000, satellite images.
\end{keyword}

\end{frontmatter}


\section{Introduction}
The new concept of natural numerical networks is introduced in this paper. It is used as a novel tool for automated classification of protected habitats using satellite images. The natural numerical network is based on the numerical solution of the nonlinear forward-backward diffusion (FBD) equation on complete graph. The {\it network}, represented by the {\em numerical} discretization of the FBD equation, can be classified as a new deep learning method. Usually, a deep learning method is formed by an artificial neural network with many hidden layers. In our discretization scheme, we use a sequence of time steps resolving the dynamics of the diffusion equation, and one hidden layer corresponds to one time step of the numerical scheme for solving the FBD equation. The proposed deep learning method does not use artificial neural network principles in its construction, see e.g. \cite{1990JGCD}, \cite{BCC} or \cite{Goodfellow-et-al-2016}. In addition, diffusion equations are widely used in modelling phenomena in natural sciences, such as biology, physics or chemistry; thus, we call the proposed network {\it natural}. Indeed, the method seems to truly be a natural procedure for the clustering and supervised classification of data, so we call the proposed method the {\it natural numerical network} or the {\it natural network} for short. In building the natural network, we are inspired by the work \cite{HR} in which Eldad Haber and Lars Ruthotto showed the relation between a successful deep learning model, the so-called Residual Neural Network (ResNet) \cite{RN, ReversArch}, and the numerical solution of the system of ordinary differential equations using the forward Euler method. Subsequently, they designed parabolic and hyperbolic networks for deep learning based on the appropriate partial differential equations \cite{HR}. Our natural network uses another type of PDEs, the nonlinear forward-backward diffusion equations. The forward diffusion causes the movement of points in a feature space toward each other. The opposite effect, when the points are kept away from each other, is caused by backward diffusion. This yields the desired classification. The natural network based on FBD equations contains just a few parameters that are optimized in the learning phase of the method. After learning the parameters and optimizing the topology of the graph of the network, the classification is performed. The relevancy map for each habitat is created as a novel tool for validating the classification, studying the relation of the habitat classification with the species composition and finding new Natura 2000 habitat appearances. It shows abilities of developed method in distinguishing between next-standing mixed deciduous forests with similar species composition and between two types of riparian forests as shown in the alluvium of the Danube river, and also potential in automated finding new localities of protected habitat areas as shown on softwood floodplain forest newly discovered in Slovakia by the proposed method.

The use of satellite data has become one of the essential methods for effectively and directly acquiring information on the Earth's surface \cite{Liu2015, Randin_etal.2020}. Together with standardized botanical records (plots) and regular in situ measurements, remote sensing is a powerful monitoring instrument \cite{Corbane_etal.2015, Lausch_etal.2020} playing an irreplaceable role in acquiring data essential for evaluating and implementing environmental policy by data analysis \cite{Liu2015, Ullerud_etal.2018, Chi_etal.2016}. Remote sensing is also one of the most important tools in ecology and nature conservation for achieving the effective monitoring of ecosystems in space and time \cite{Rocchini}; thus, using satellite images to monitor habitats and biome dynamics has been highlighted in many types of research activities. Ecosystems representing Natura 2000 habitats are complex plant communities including tree, shrub, and herb layers together with typical fauna \cite{EEANatura2000, Natura2000} for which it was impossible to reach an automated identification and classification with existing methodologies based on satellite data \cite{Bock, Vanden} so far.
In this paper we develop method which performs the first step to fill this gap. Although the method is generally designed to work with any type of optical data monitoring the Natura 2000 habitats, we use the optical information from spectral bands of the Sentinel-2 satellite \cite{ESASentinel} freely available on European Space Agency (ESA) servers.

In summary, the aims of presented study are 
1) to give the complete mathematical definition of the natural numerical network based on forward-backward diffusion equations, 
2) to present the numerical scheme corresponding to the natural numerical network in all phases of the classification algorithm,
3) to optimize the natural numerical network on training dataset of Natura 2000 habitats, and 
4) to apply the trained natural numerical network to classification of Natura 2000 habitats, construction of the habitat relevancy maps for results validation and for finding new appearances of protected habitats.

\section{Methods}
\subsection{Mathematical model}
	Let us define a graph as an ordered pair $ G = (V(G), E(G)) $, where $ V(G) $ is a finite set of vertices, and $ E(G) $ is a set of the two-element subsets of $ V(G) $ representing the edges of the graph $ G $ \cite{teoriaGrafovKnor}. We denote the number of vertices of the graph $ G $ by $ N_V $. Let us consider that the graph $ G $ is a complete graph that means that each vertex $ v \in V (G) $ is connected to each other vertex by an edge. Let us suppose that graph $ G $ is undirected and thus the edges do not have an orientation.

	Let us consider the function $ X : G \times [0, T] \rightarrow R^k $ representing the spatial coordinates $ X(v, t) = (x_1(v, t), \ldots, x_k(v, t)) $ of the vertex $ v \in V(G) $ in time $ t \in [0, T] $. In our case, $ k $ is a dimension of feature space $ R^k $. A diffusion of $ X(v, t) $ on the graph $ G $ is formulated as a partial differential equation (PDE)
\begin{equation} 
\label{eq:nonlinearPDE}
	\partial_{t}X(v, t) = \nabla \cdot (g \nabla X(v, t)), \qquad v\in V(G), \qquad t\in [0, T],
\end{equation}
where $ g $ represents a so-called diffusion coefficient, see also \cite{calculusOnGraphs}. We consider equation (\ref{eq:nonlinearPDE}) together with an initial condition $ X(v, 0) = X^0(v) $, $ v \in V(G) $. The boundary conditions are not necessary to prescribe because diffusion occurs between all vertices of the complete undirected graph $ G $.

	Let us define the distance between two vertices $ v, u \in V(G) $ as the Euclidean distance between two points $ X(v, \cdot) $ and $ X(u, \cdot) $ in the feature space $ R^k $ and denote it $ L(v, u) $. Since every pair of vertices of $ G $, $ v $ and $ u $, forms one edge, $ e = \{v, u\} $, we can simplify the notation for the distance and denote it as a length of the edge $ L(e) $, which will be used throughout the following text.

\begin{figure}
 \begin{center}
	\includegraphics[width= 0.6 \textwidth]{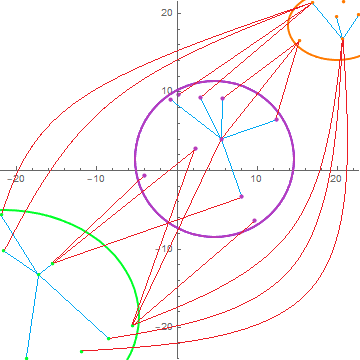}
	\caption{Randomly generated 2D points in three clusters and some links of forward diffusion (blue lines) inside the clusters, and some links of backward diffusion (red lines) between points from different clusters.}
	\label{fig:dopredSpat-linky}
 \end{center}
\end{figure}
	
	We will design the diffusion coefficient $ g $ depending on the length of the edges  $ e $ of the graph $ G $. It will give a nonlinear diffusion model on the graph representing a generalization of the Perona-Malik model from the image processing \cite{Perona-Malik}. We consider equation (\ref{eq:nonlinearPDE}) with the diffusion coefficient $ g $ in the form
\begin{equation} 
\label{eq:difKoef}
	g(e) = \varepsilon(e) \frac{1}{1 + K L^2(e)}, \qquad K \geq 0.
\end{equation}
The value of  $ \varepsilon(e) $ in the diffusion coefficient depends on the type of diffusion which is applied on every single edge. If we need to apply forward diffusion, we choose $ \varepsilon(e) $ as a positive constant. The forward diffusion, represented by a positive diffusion coefficient, averages the values of a diffused quantity. The averaging reflects the smoothing property of the standard diffusion equation. It is called forward because it describes diffusion towards the future and in our application, it causes a moving and thus clustering of points together. On the other hand, the backward diffusion is represented by a negative diffusion coefficient, $ \varepsilon(e) $ is a small negative value, and can be understood as returning to the past in a diffusion process. It is an inverse process of the smoothing (averaging) values of a diffused quantity and in our application, it gives a repulsion of the points belonging to different clusters. If we only used the backward diffusion model, the points would be moving away from each other, and the whole system became unstable. But by a suitable combination of the forward and backward diffusion, when we choose a small negative coefficient $ \varepsilon(e) $ for backward diffusion, we do not observe any calculation instability. Such a model is a suitable and natural tool for supervised learning - the points inside the given clusters are kept together while points of different clusters are kept away from each other. Such proper behaviour is realized by the model (\ref{eq:nonlinearPDE})-(\ref{eq:difKoef}).

	In Fig. \ref{fig:dopredSpat-linky} we illustrate behaviour of the model (\ref{eq:nonlinearPDE})-(\ref{eq:difKoef}) on three given clusters where by blue lines we plot some of the links of forward diffusion and by red lines some of the links of backward diffusion. This figure depicts the basic features and behaviour of the natural network. The points inside a given cluster are attracted by the forward diffusion while there is a repulsion of points of different clusters by the backward diffusion.
	
\begin{figure}
 \begin{center}
	\includegraphics[width= 0.6 \textwidth]{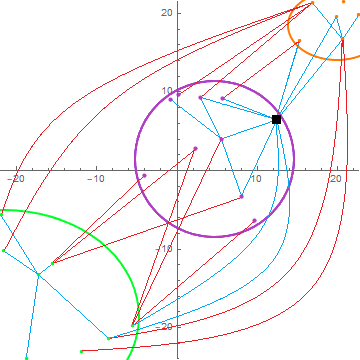}
	\caption{Randomly generated 2D points in three clusters with one  new observation (black square). The blue lines represent some of the forward diffusion links inside the clusters and between the new observation and all other points. The red lines represent the links of backward diffusion between points from different clusters.}
	\label{fig:dopredSpat-linky_newcom}
 \end{center}
\end{figure}

	Additionally, in Fig. \ref{fig:dopredSpat-linky_newcom} we illustrate the situation arising in supervised learning and application phases of the classification method when a new observation is added into the network. Only the forward diffusion is applied to all links of the vertex representing the new observation. It is depicted by the blue lines connecting the new observation (black square) with every other point. Thus, this new observation is attracted by a certain diffusion speed to all existing clusters which themselves are subject to the forward-backward diffusion as described before. The dynamics of the network decides about the cluster membership of the new observation.

	The model (\ref{eq:nonlinearPDE})-(\ref{eq:difKoef}) contains, together with forward-backward diffusion switch $ \varepsilon(e) $, also a weighting coefficient $ K $. The constant $ K $ controls how the length $ L(e) $ of the edge $ e = \{v, u \} $ affects the diffusion of the vertices $ v $ and $ u $ over time. If $ K L(e)^2 $ is large, the diffusion coefficient $ g $ is close to 0 which means that the diffusion process will be slow, and the points are not diffusing (moving to each other) by an averaging. If $ K L(e)^2 $ is small, the diffusion coefficient is close to 1, the diffusion process is faster and points are moving to each other fast by the diffusion. Since the coefficient $ K $ is multiplying the squared length of the edge, the distant points in the feature space are averaging (moving) slower than the close points. 
	
	The model (\ref{eq:nonlinearPDE})-(\ref{eq:difKoef}), after discretization, represents a basic natural numerical network for supervised deep learning classification and, with just a positive coefficient $ \varepsilon(e) $, it is also a proper model for unsupervised clustering (which is, however, not discussed in this paper). This basic model allows useful modifications which will be used also in our final classification algorithm. First of all, we slightly modify the diffusion coefficient into the form
\begin{equation} 
\label{eq:difKoefExtend}
	g(e) = \varepsilon(e) \frac{1}{1 + \sum_{i = 1}^{k} (K_i \ l_i^2(e))}, \qquad K_i \geq 0.
\end{equation}
Now, the parameters $ K_i $, $ i = 1, \ldots, k $, represent weights for each coordinate $ l_i(e) $, $ i = 1, \ldots, k $, of the vector 
\begin{equation}
\begin{aligned}
	l(e) &= (l_1(e), \ldots, l_k(e))^T = X(u, \cdot) - X(v, \cdot) = \\ &= (x_1(u, \cdot) - x_1(v, \cdot), \ldots, x_k(u, \cdot) - x_k(v, \cdot))^T, \qquad v, u \in V(G).
\end{aligned}
\end{equation} 
By this modification, we can control the diffusion speed in each direction of the $ k $-dimensional feature space and achieve more accurate classification results.

\begin{figure}
 \begin{center}
	\includegraphics[width= 0.6 \textwidth]{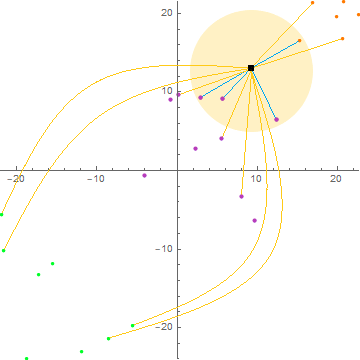}
	\caption{Randomly generated 2D points in three clusters with one new observation (black square). A "$ \delta $-diffusion neighborhood" around the new observation (yellow circle) and some links of forward diffusion with non-zero value of diffusion coefficient $ g $(blue lines), and some links where the diffusion is set to zero value (yellow lines), because the points are outside of the $ \delta$-diffusion neighborhood.}
	\label{fig:d-okolie}
 \end{center}
\end{figure}

	Next modification of the basic model (\ref{eq:nonlinearPDE})-(\ref{eq:difKoef}) controls a forward diffusion coefficient on the edges of the new observation points, see Fig \ref{fig:dopredSpat-linky_newcom} and \ref{fig:d-okolie}. We can reduce the forward diffusion influence on the vertex $ v \in V(G) $ by using the diffusion coefficient in the form
\begin{equation} 
\label{eq:difKoefDelta}
	g(e) = \max(\varepsilon(e) \frac{1}{1 + \sum_{i = 1}^{k} (K_i \ l_i^2(e))} - \delta, \ 0), \qquad \varepsilon(e) > 0 
\end{equation}
on all edges $e$ of $ v $, where $ \delta $ is a parameter of the "diffusion neighborhood" size.  In classification of the new observation, the aforementioned modification causes that only the points in a "$ \delta $-diffusion neighborhood", i.e. for which the diffusion coefficient is large than $\delta$, are attracting new observation point $ v $. This modification is illustrated in Fig. \ref{fig:d-okolie}.

\subsection{Numerical discretization - natural network construction}
	Let us denote by $ f(v, t) $ any of the coordinates $ x_i(v, t) $ of $ X(v, t) = (x_1(v. t), \ldots, x_k(v, t)) $. To discretize the equation (\ref{eq:nonlinearPDE}), we use {\em i)} the balance of diffusion fluxes (inflows and outflows) in each vertex $ v \in V(G) $ and {\em ii)} the approximation of the diffusion flux to the vertex $ v $ along its edge $ e $.
	
	First, let us define the diffusion flux approximation, which depends on the difference of values of function $ f $ at the vertices $ v $ and $ u $, as
\begin{equation} 
\label{eq:aproxDifTok}
	\mathcal{F}(v, e, t) = g_e (f(u, t) - f(v, t)),
\end{equation}
for each edge $ e = \{v, u\} $, where $ g_e $ represents the diffusion coefficient on the edge $ e $. If $ \mathcal{F}(v, e, t) > 0 $, it represents the diffusion inflow of the quantity $ f $ into the vertex $ v $. On the other hand, if $ \mathcal{F}(v, e, t) < 0 $, it represents the diffusion outflow of the quantity $ f $ from the vertex $ v $. Then the balance of diffusion fluxes in vertex $ v $ is expressed by the equation
\begin{equation} 
\label{eq:bilanciaDifTok}
	\partial_{t} f(v, t) = \sum_{e \ni v} \mathcal{F}(v, e, t),
\end{equation}
that means, the time derivative of $ f $ in the vertex $ v $ is positive - the value of $ f $ increases in time, if the overall inflow to the vertex $ v $ is greater than the overall outflow from the vertex $ v $. Vice versa, the time derivative of $ f $ in the vertex $ v $ is negative if the sum of inflows and outflows in the vertex $ v $ is negative, i.e. outflows from the vertex $ v $ are greater than inflows. When we substitute the approximation of the diffusion flux (\ref{eq:aproxDifTok}) into the balance equation (\ref{eq:bilanciaDifTok}), we obtain
\begin{equation} 
\label{eq:bilanciaAproxDifTok}
	\partial_{t} f(v, t) = \sum_{\substack{e \ni v \\ e = \{v, u\}}} g_e (f(u, t) - f(v, t)).
\end{equation}
The right hand side of the equation (\ref{eq:bilanciaAproxDifTok}) in the graph theory represents the so-called "graph-Laplacian" (see equation (12) in \cite{calculusOnGraphs} or equation (2.5) in \cite{edge-basedLaplacian}) which is given for a weighted complete undirected graph by relation
\begin{equation} 
\label{eq:graphLaplacian}
	\nabla \cdot (\nabla f)(v, t) = \frac{1}{\nu(v, t)} \sum_{\substack{e \ni v \\ e = \{v, u\}}} g_e (f(u, t) - f(v, t))
\end{equation}
where $ \nu(v, t) $ represents a measure of the vertex $ v $ and $ g_e $ represents a "weight" of the edge in the weighted graph. In numerical mathematics, we would understand the "Laplacian" defined in this way as an averaged Laplace operator on a "finite volume" $ v $ with a  measure (area/volume) $ \nu (v, t) $. Our numerical discretization (\ref{eq:bilanciaAproxDifTok}) of the diffusion equation on the complete undirected graph corresponds to the choice $ \nu(v, t) = 1 $, which is the standard choice for the vertex measure also in the graph theory, see \cite{calculusOnGraphs}.

\begin{figure}
 \begin{center}
	\includegraphics[width = 0.4\linewidth]{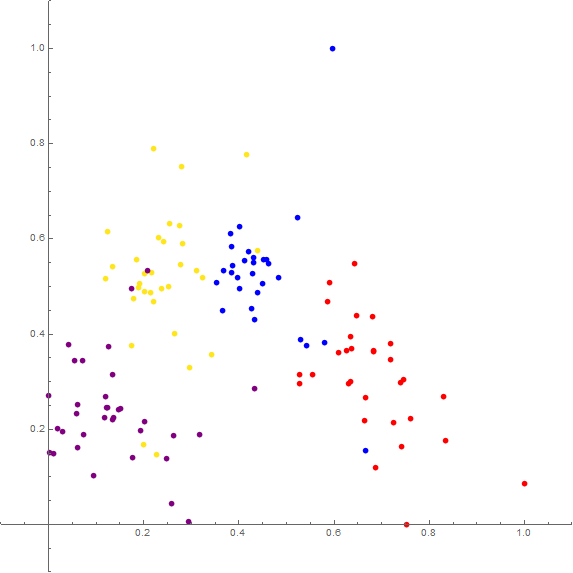}~~
	\includegraphics[width = 0.4\linewidth]{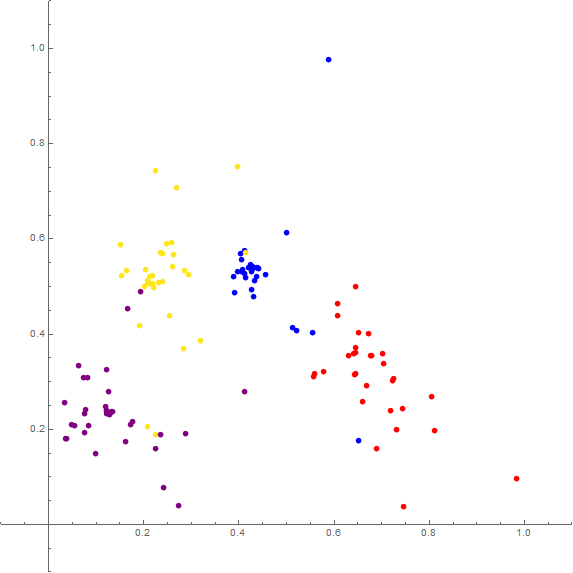} \\
	\includegraphics[width = 0.4\linewidth]{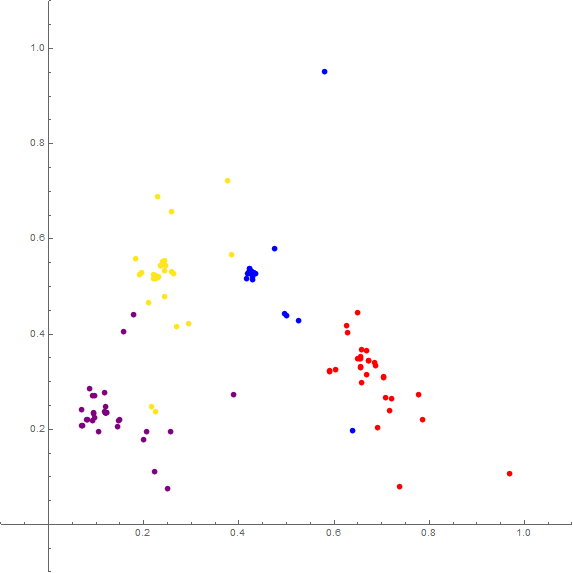}~~
	\includegraphics[width = 0.4\linewidth]{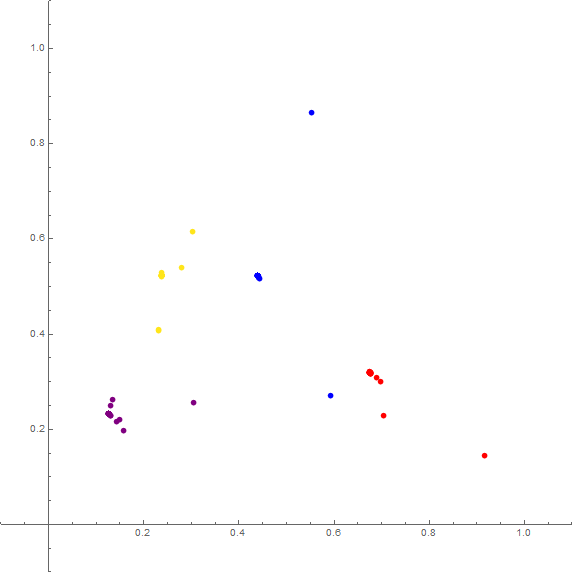} \\
\hskip 1mm	\includegraphics[width = 0.4\linewidth]{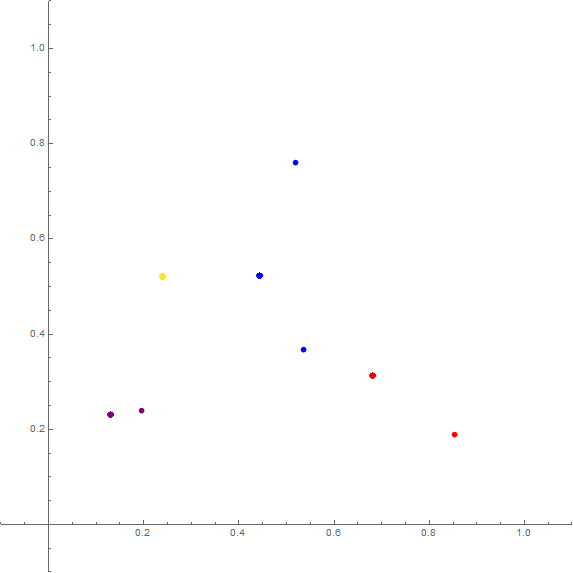}~~
	\includegraphics[width = 0.4\linewidth]{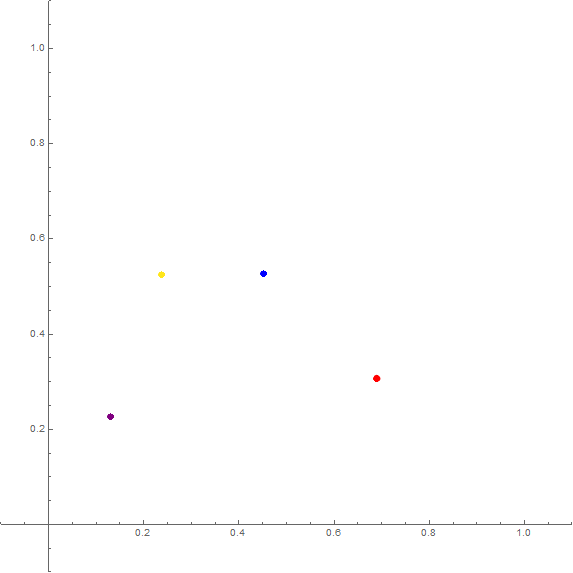} 
	\caption{A dataset with $ 125 $ observations to which the forward-backward diffusion natural network is applied. We show the dynamics of the network in time steps $ n = 0, 1, 2, 5, 8, 14 $. The observations are classified into $ 4 $ clusters and the points in different clusters are distinguished by different colors.}
	\label{fig:125x72_Squares_rev1}
 \end{center}
\end{figure}

\begin{figure} 
 \begin{center}
	\includegraphics[width = 0.4\linewidth]{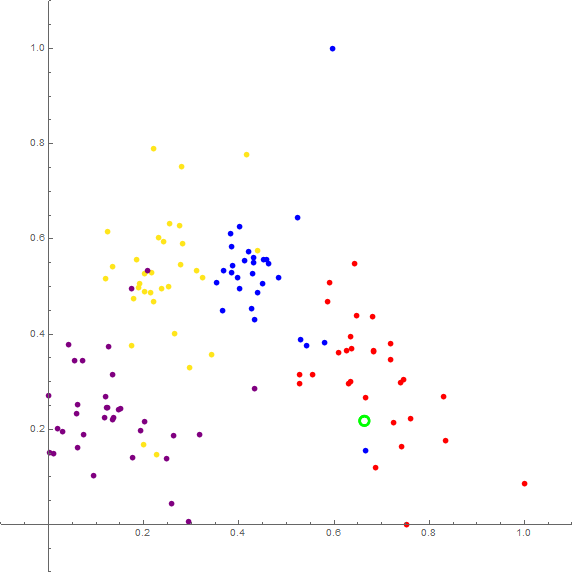}~~
	\includegraphics[width = 0.4\linewidth]{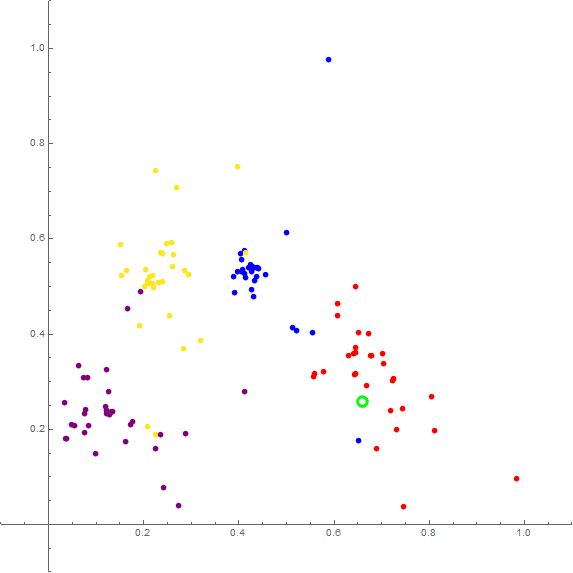} \\
	\includegraphics[width = 0.4\linewidth]{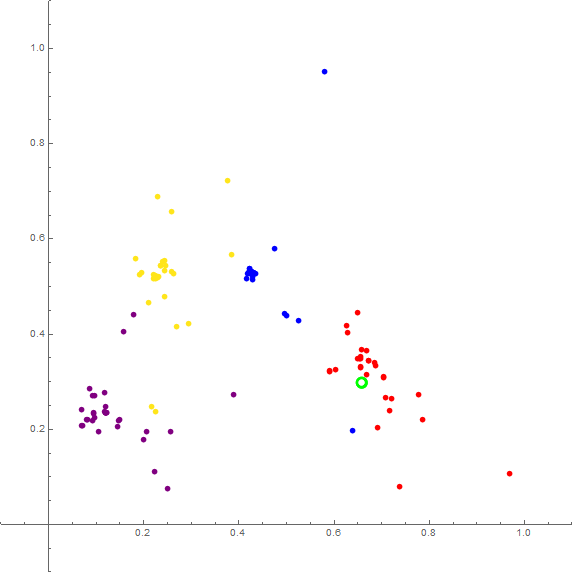}~~
	\includegraphics[width = 0.4\linewidth]{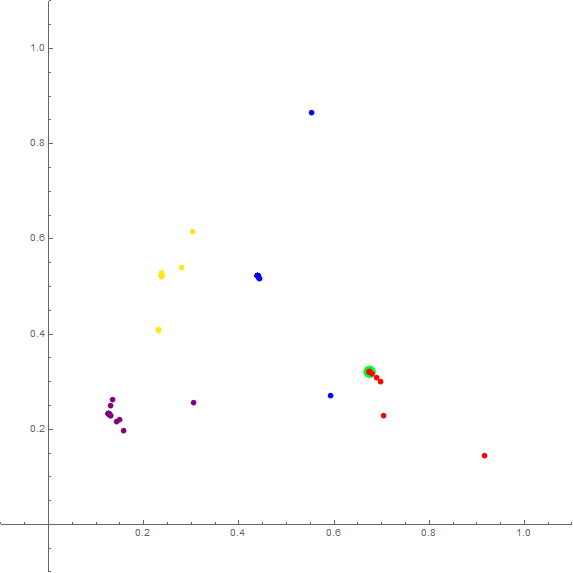} \\
	\hskip 1mm\includegraphics[width = 0.4\linewidth]{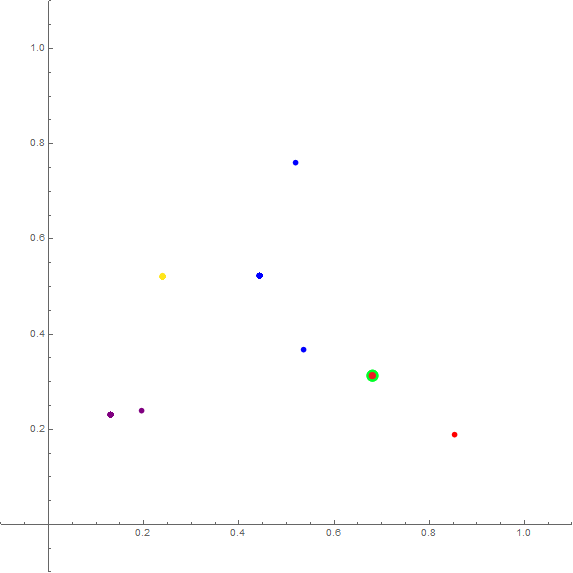}~~
	\includegraphics[width = 0.4\linewidth]{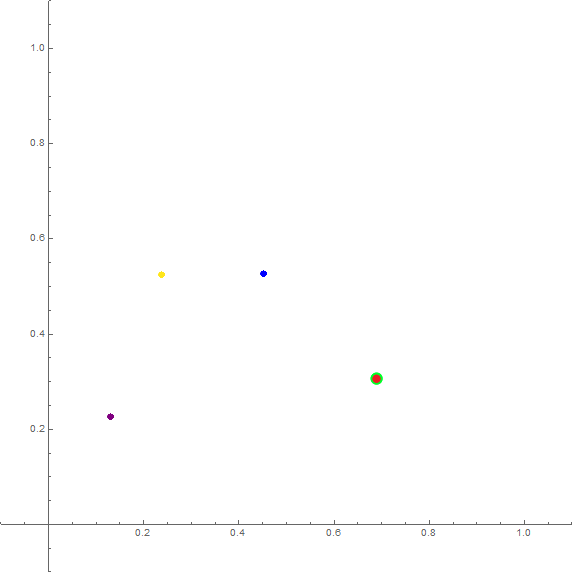} 
	\caption{A dataset with $ 124 $ observations and one "new" observation (green point) to which the forward-backward diffusion is applied. We show the dynamics of the network in time steps $ n = 0, 1, 2, 5, 8, 14 $. The new observation is classified to the "red" cluster. This figure illustrate the learning phase of the classification.}
	\label{fig:125x72_Squares_rev1_n0}
 \end{center}
\end{figure}

	For the time discretization, we use the semi-implicit approach, see e.g. \cite{semiImpl}. The time interval $ [0, T] $ is divided uniformly into $ M $ time steps $ t_{n}, \ n = 1, \ldots, M $ and let $ \tau $ denote the size of the time step. For the approximation of time derivative we use the finite difference method and obtain 
\begin{equation} 		
\partial_t f(v, t) = \frac{f^n(v) - f^{n - 1}(v)}{\tau},
\end{equation}
where $ f^n(v) = f(v, t_n) $.
	
	Since the diffusion coefficient $ g_e $ at the edge $ e = \{v, u\} $ can depend on the unknown quantity $ f $, see (\ref{eq:difKoef})-(\ref{eq:difKoefDelta}), and thus can change over time, we take its value from the previous time step. We denote it by $ g^{n - 1}_e $ and obtain the semi-implicit scheme in the form
\begin{equation}
\label{eq:poAplSpatDif}
\frac{f^n(v) - f^{n - 1}(v)}{\tau} = \sum_{\substack{e \ni v \\ e = \{v, u\}}} g^{n - 1}_e (f^n(u) - f^n(v)).
\end{equation}
The semi-implicit scheme (\ref{eq:poAplSpatDif}) can by rewritten in each time step $ n = 1, \ldots, M $ into the system of linear equations
\begin{equation}
\label{eq:sysLE}
(1 + \tau \sum_{\substack{e \ni v \\ e = \{v, u\}}} g^{n-1}_e) f^n(v) - \tau \sum_{\substack{e \ni v \\ e = \{v, u\}}} g^{n-1}_e f^n(u) = f^{n - 1}(v).
\end{equation}
This system of equations is represented by a full matrix and as we have said before, for a complete undirected graph it is not necessary to define any boundary condition. 

	In the case of classification of the data from $ k $-dimensional feature space, our diffusing variables are the Euclidean coordinates $ X(v, t) = (x_1(v, t), \ldots, x_k(v, t)) $ of the vertices $ v $ of the graph $ G $. In general, we get in each time step $ k $ systems of linear equations
\begin{equation}
\label{eq:kSysLE}
	(1 + \tau \sum_{\substack{e \ni v \\ e = \{v, u\}}} g^{n-1}_e) x^n_i(v) - \tau \sum_{\substack{e \ni v \\ e = \{v, u\}}} g^{n-1}_e x^n_i(u) = x^{n - 1}_i(v), \qquad i = 1, \ldots, k, \quad v \in V(G),
\end{equation}	
which are interconnected by the diffusion coefficient $ g^{n-1}_e $, which depends on all $ x_i^{n - 1}(v) $, $ x_i^{n - 1}(u) $, $ i = 1, \ldots, k $ and can be written in the form
\begin{equation} 
\label{eq:difKoefNum}
	g_e^{n - 1} = \varepsilon(e^{n - 1}) \frac{1}{1 + \sum_{i = 1}^{k} (K_i \ l_i^2(e^{n - 1}))}, \qquad K_i \geq 0
\end{equation}


Let us denote the $i$-th cluster by $ C_i $, $ i = 1, \ldots, N_C $, where $ N_C $ is the number of clusters. Let us have the vertices $ v \in C_l $ and $ u \in C_m $, where $ l, m \in \{1, \ldots, N_C\} $, $ e^{n - 1} = \{v, u\} $. The value $ \varepsilon(e^{n - 1}) $ in the diffusion coefficient (\ref{eq:difKoefNum}) is given by the following values

\begin{equation} 
\label{eq:epsForFBD}
 \begin{array}{ll}
	\varepsilon(e^{n - 1}) \geq 0, \qquad {\rm if} \quad l = m,  \\
	\varepsilon(e^{n - 1}) < 0, \qquad {\rm if} \quad l \neq m.
    			\end{array}
\end{equation}

	For a reader convenience in Fig. \ref{fig:125x72_Squares_rev1} we present dynamics of the network by (\ref{eq:kSysLE}) - (\ref{eq:epsForFBD}) in order to show how the points inside the given clusters are moving together and clusters themselves are keeping away. In practice, in the learning phase and also in the application phase, we add new observation to the network and run the dynamics of the network. In the application phase we add to the network completely new observation while in the learning phase the new observation is taken away from the learning dataset. In both cases, the dynamics is modified in such a way that all other points are moving by (\ref{eq:kSysLE}) - (\ref{eq:epsForFBD}) but for the new observation $ w \notin C_i $, $ i \in \{1, \ldots, N_C\} $, $ e = \{w, u\} $, the diffusion coefficient is set to 
\begin{equation} 
\label{eq:difKoefNum_newcom}
	g_e^{n - 1} = \max(\varepsilon(e^{n - 1}) \frac{1}{1 + \sum_{i = 1}^{k} (K_i \ l_i^2(e^{n - 1}))} - \delta, \ 0), \quad \varepsilon(e^{n - 1}) \geq 0, \ K_i \geq 0, \ \delta > 0
\end{equation}
are given constants.

\begin{figure}[t] 
 \begin{center}
	\includegraphics[width = 0.6\linewidth]{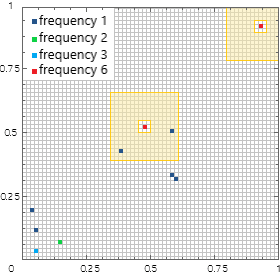}
	\caption{For illustration of the histogram stopping criterion we visualize the 2D grid with spacing $h = 0.015625 $ (grey lines), the marked grid cells (red) in which frequency is greater or equal to $ S_{min} = 6 $, the neighborhood of the marked cell which is also examined (yellow subdomains) and further cells coloured by their frequencies given in the upper left corner. We see that the cluster was already formed inside the marked cell in the upper right corner because there are no other points inside the examined neighborhood while the cluster is not yet formed in the central marked cell because there are still other points in the examined neighborhood, there are two cells with the frequency equal to $ 1 $.}
	\label{fig:histogram2D}
 \end{center}
\end{figure}

	A stopping criterion is applied for dynamics of the network. We call it the histogram stopping criterion because it calculates the number of occurrences (frequency) of evolving points in prescribed spatial cells in every time step. That means the $ k $D grid with cells given by a specific spacing $h$ is created, see the grey grid in Fig. \ref{fig:histogram2D} for 2D case illustration. In numerical experiments presented below, we always use the spacing $h=0.01$. We also determine which of the given clusters is the smallest one, and let the variable $ S_{min} $ represent the number of points in this smallest cluster. Then, in every time step the histogram development is monitored, and whenever the number of points (denoted as the frequency in Fig. \ref{fig:histogram2D} left up corner) inside a grid cell is greater or equal to $ S_{min} $, the cell is marked, see the red cells in Fig. \ref{fig:histogram2D}. At the same time, a specific neighborhood of every marked cell is examined, see the yellow subdomains in Fig. \ref{fig:histogram2D}. The examined neigborhood is given by the interior of the concentric squares with Chebyshev radius $ H_1 $ and $ H_2 $, respectively, in Fig. \ref{fig:histogram2D} we illustrate the situation where $H_1=1$ and $H_2=8$, and such parameters are also used in the computations presented below. If there are only zero values of the frequency in all the cells inside the examined neighborhood of the marked cell, we claim that a cluster was formed in the marked cell. The dynamics of the network is stopped when the number of clusters formed is equal to the number of clusters given. 	

\begin{figure}[t] 
 \begin{center}
	\includegraphics[width = 0.48\linewidth]{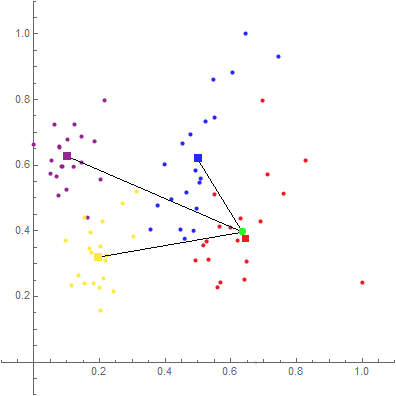}~~
	\includegraphics[width = 0.48\linewidth]{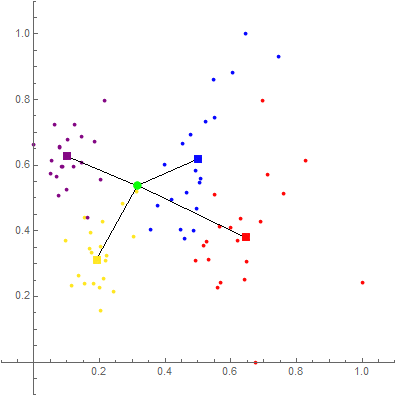}
	\caption{We show data in the time step $ n = 0 $, with four given clusters (various colors) and one new observation (green point), the centroids of the given clusters (colored squares) and distances of the new observation to the centroids of the given clusters (black lines). In the left picture, the new observation resulting classification relevancy will be high, close to 1, while in the right picture it will be low, around 0.5.}
	\label{fig:DS80newcom1a20rel1}
 \end{center}
\end{figure}

	For quantifying the relevancy of classification of any new observation $ w \notin C_i $, $ i \in \{1, \ldots, N_C\} $, the relevancy coefficient $ R(w)$ is defined. We define it by using the information to which cluster the new observation is classified combined with the information about distances of the point representing the new observation and the centroids of all given clusters. 
	
The new observation is classified into the cluster $C_a$, if after stopping the network dynamics its closest point in the network is a point from the cluster $C_a$ and, at the same time, their distance is less than $H=10 h$. Otherwise, the point is not classified to any cluster, it is called the outlier and its relevancy will be equal to 0 in all clusters.
	 
\begin{figure}[H] 
 \begin{center}
	\includegraphics[width = 0.5\linewidth]{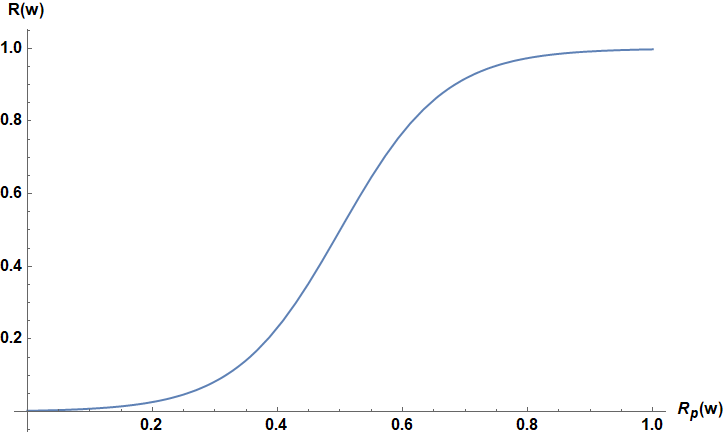}
	\caption{The relevancy coefficient $R(w)$ plotted in dependence on $R_p(w)$, $\lambda=12$ in the definition of the logistic function 
	{\rm (\ref{eq:logistic})}.}
	\label{fig:logisticRelevancy}
 \end{center}
\end{figure}

Now, let us assume that the new observation is classified into the cluster $C_a$. The new observation has given its initial position $ X(w, 0) $ and we calculate the centroids of the given clusters at the final time step by
\begin{equation}
	{\mathcal{C}}_i = \frac{1}{N_{C_i}} \sum_{v \in C_i} X(v, T), \qquad i = 1, \ldots, N_C,
\end{equation}
where $ N_{C_i} $ is the number of points in the cluster $ C_i $. Then we calculate the distance between the new observation and the centroid of the cluster $ {\mathcal{C}}_a = {\mathcal{C}}_a(w) $ to which it is assigned by the network dynamics,
\begin{equation}
	l_1(w) = \mid X(w, 0) - {\mathcal{C}}_a\mid\ ,
\end{equation}
and the average distance of the new observation to all other cluster centroids,
\begin{equation}
	l_2(w) = \frac{1}{N_c - 1} \sum_{\substack{i=1 \\ i \neq a}}^{N_c} \mid X(w, 0) - {\mathcal{C}}_i \mid\ .
\end{equation}
The above distances are used to define the quantity
\begin{equation}
	\label{eq:relevancy}
	R_p(w) = 1 - \frac{l_1(w)}{l_1(w) + l_2(w)}
\end{equation}	
which is in the range $ [0, 1] $ and which is the basis for definition of the relevancy coefficient. When the position of the new observation $ X(w, 0) $ is close to the centroid of the cluster to which it is assigned then $R_p(w)$ is close to 1, see Fig. \ref{fig:DS80newcom1a20rel1} left. In this case, the relevancy of the resulting classification should be high. The quantity $R_p(w)$ is close to $0.5$ or less if the distance of the new observation to the centroid of the cluster to which it is assigned is similar or greater than its distance to other clusters centroid, see Fig. \ref{fig:DS80newcom1a20rel1} right. In this case, the resulting relevancy of the classification should be significantly reduced. We use the logistic function 
\begin{equation}
	\label{eq:logistic}
    \mathcal{L} (x) = \frac{1}{1+e^{\lambda (0.5-x)}}
\end{equation}	
which after linear rescaling from the interval $[\mathcal{L}(0), \mathcal{L}(1)]$ to the interval $[0,1]$ give the final definition of the relevancy coefficient $R(w)$ for any new observation $w$,
\begin{equation}
	\label{eq:relevancy1}
    R(w) = \frac{\mathcal{L}(R_p(w)) - \mathcal{L}(0)}{\mathcal{L}(1)- \mathcal{L}(0)}.
\end{equation}	
While $R_p(w)$ is linearly decreasing from 1 to 0, depending on the distances ratio $l_1(w)/(l_1(w) + l_2(w))$, the final relevancy coefficient $R(w)$ has the nonlinear character, see Fig. \ref{fig:logisticRelevancy}. It sets the relevancy values to be close to 1 for all new observations belonging to a neighborhood (size of which depends on $\lambda$) of the centroid of the cluster to which it is assigned. In all numerical experiments presented below, we use $\lambda=12$. The relevancy coefficient defined by (\ref{eq:relevancy1}) is used in the definition of relevancy maps which give a useful information on Sentinel-2 image pixels and image subareas membership in respective clusters, see Tabs. \ref{tab:meanRCE0} - \ref{tab:meanRC10}.

\subsection{Ground-based vegetation data sampling - Application to habitat identification and prediction}	
	
One of the possible applications of the natural numerical networks described above is ecology and nature conservation. The identification and classification of the Natura 2000 habitats by using remote sensing data still has strong limitations due to habitats’ natural character and variability. In our application, we use a natural numerical network for the classification of the protected Natura 2000 forest habitats in the territory of Slovakia. There were four Natura 2000 forest habitats dominant in Western Slovakia chosen for the classification, as shown in Table \ref{tab:habitats}. Therefore, we have four clusters $C_i$, where $i=1,\dots,N_C$, and $N_C=4$. All habitat areas borders were semi-automatically segmented in NaturaSat software \cite{segSemiAuto,Ambroz} and checked in the field by botany experts during vegetation seasons of 2019 and 2020.

\begin{figure}
 \begin{center}
	\includegraphics[width = 0.9\linewidth]{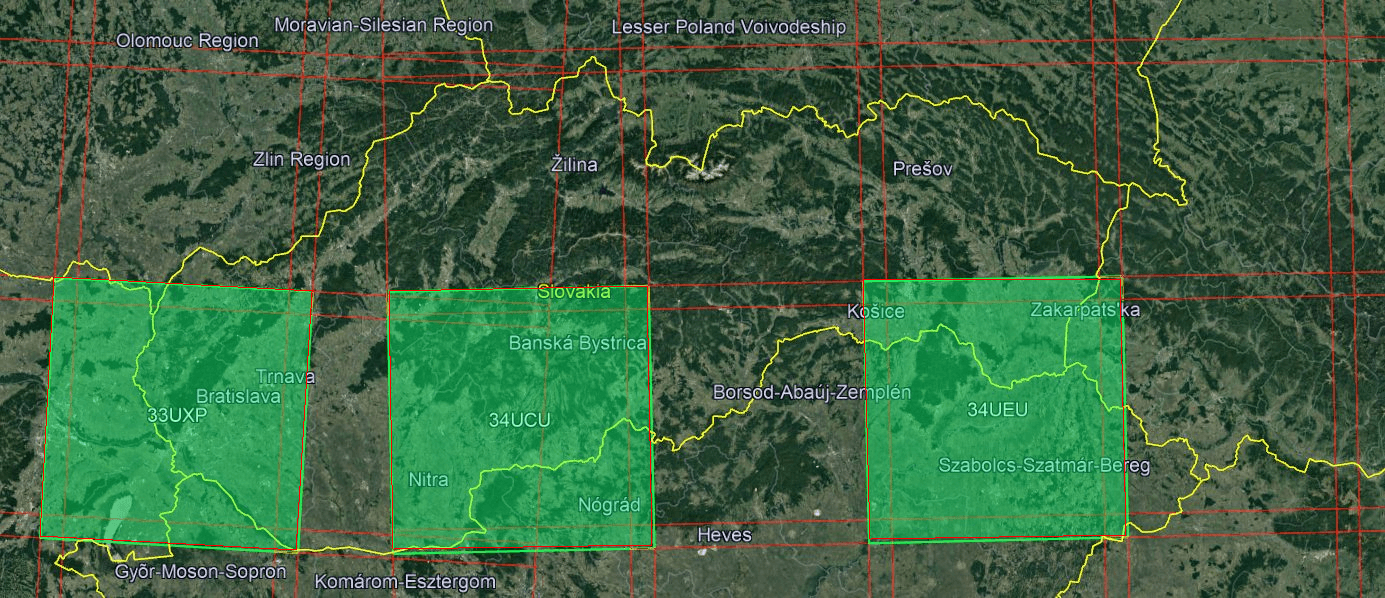}
	\caption{The map with the squares covered by the Sentinel-2 data. The left light green square {\rm Sentinel-2 tile number 34UXP} covers Western Slovakia, the middle square {\rm Sentinel-2 tile number 34UCU} covers the south of Central Slovakia and the right square {\rm Sentinel-2 tile number 34UEU} covers the south of Eastern Slovakia. 
	}
	\label{fig:SK}
 \end{center}
\end{figure}

\begin{figure}
 \begin{center}
	\includegraphics[width = 0.48\linewidth]{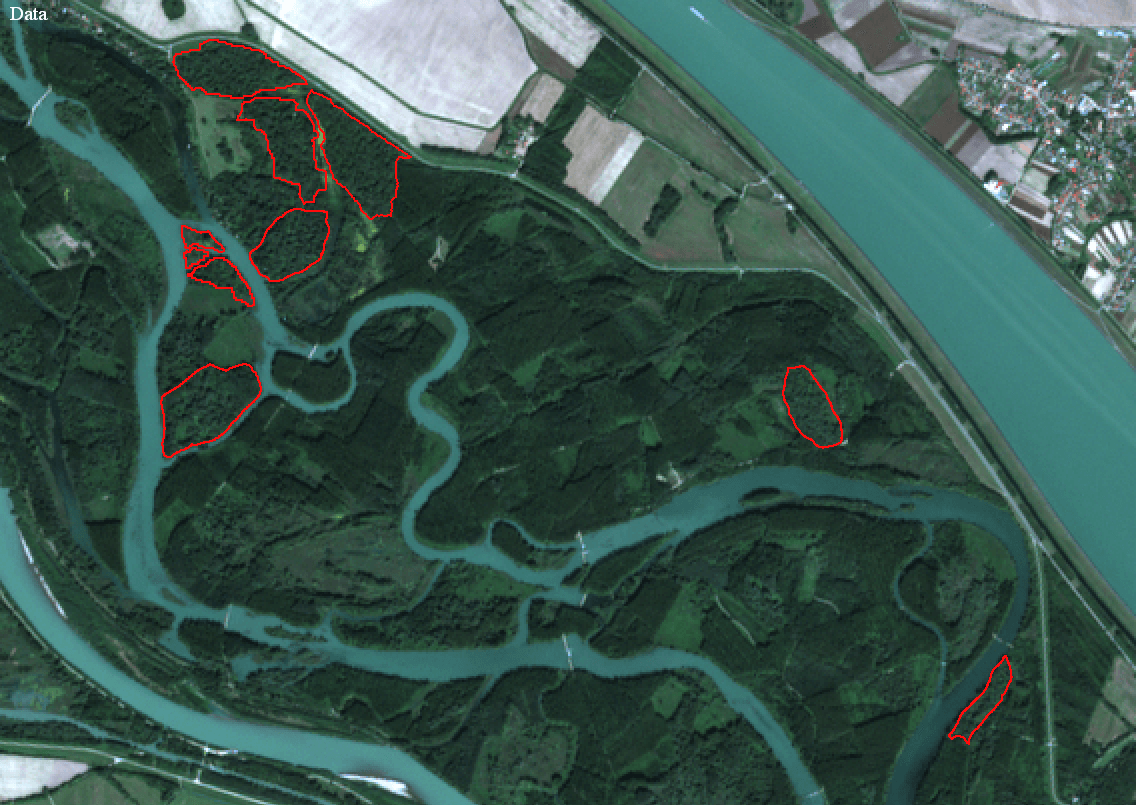}~~
	\includegraphics[width = 0.48\linewidth]{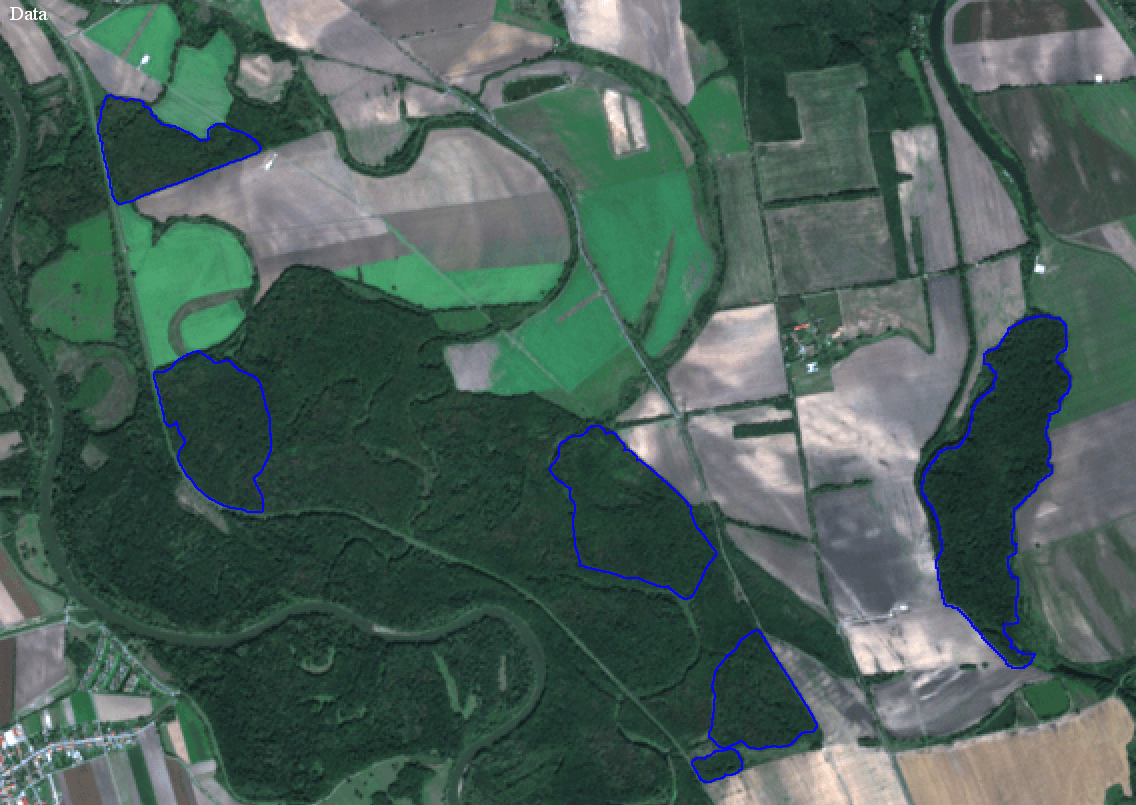} \\
	~\includegraphics[width = 0.48\linewidth]{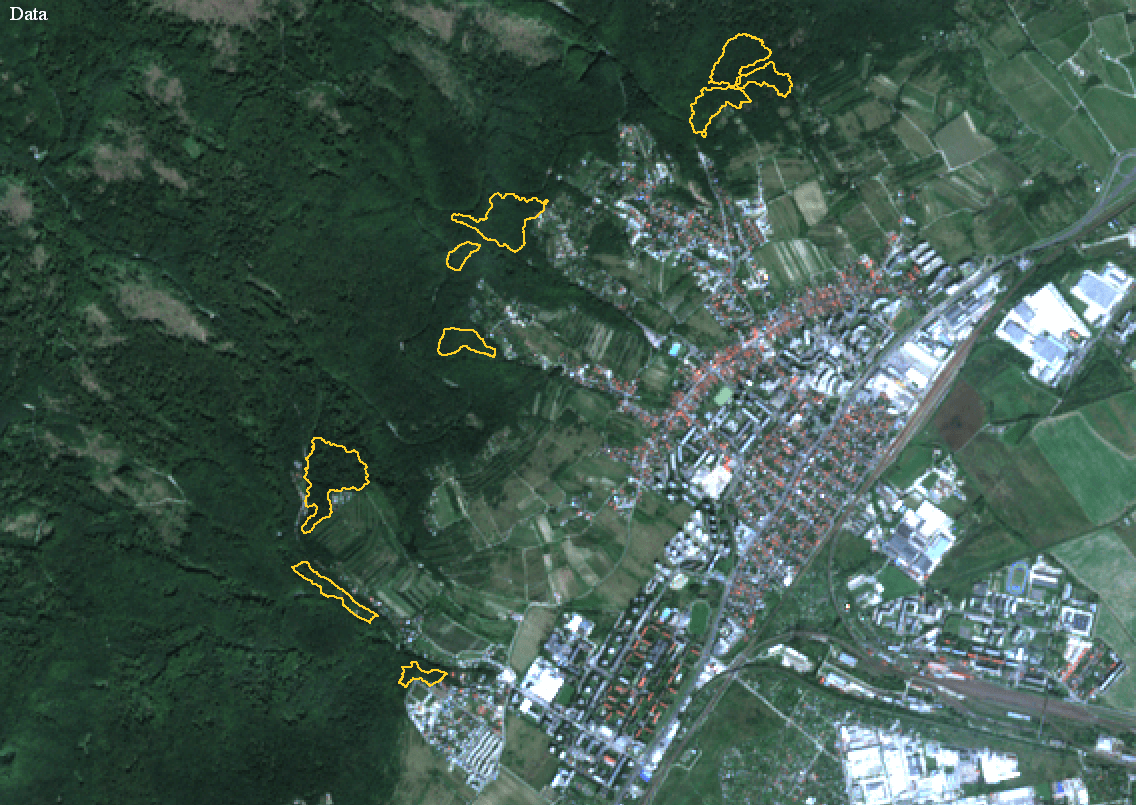}~~
	\includegraphics[width = 0.48\linewidth]{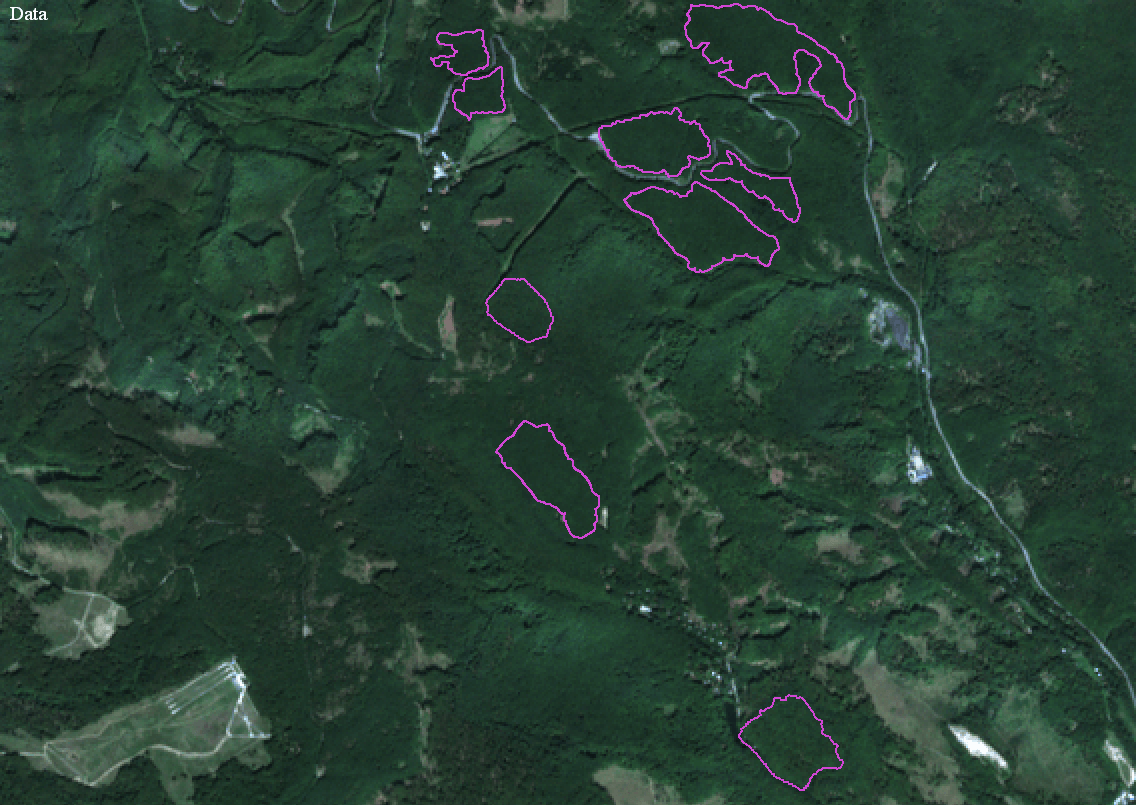}
	\caption{The subregions of Western Slovakia with segmented areas of protected Natura 2000 habitats. The segmented areas of {\rm 91E0} habitat are plotted by red curves, {\rm 91F0} habitat by blue curves, {\rm 91G0} habitat by yellow curves and {\rm 9110} habitat by purple curves.}
	\label{fig:habitats}
 \end{center}
\end{figure}

\begin{table}[H]
\begin{center} 
\begin{tabular}{ c c c }
	\hline
	Code & Name of habitat & Color \\
	\hline
	\hline
	91E0 & Alluvial forests with Alnus glutinosa and Fraxinus excelsior & red \\
	\hline 
	91F0 & \begin{tabular}{@{}c@{}} 
			Riparian mixed forests of Quercus robur, Ulmus laevis and \\
			Ulmus minor, Fraxinus excelsior or Fraxinus angustifolia, \\
			along the great rivers
		   \end{tabular} & blue \\
	\hline 
	91G0 & Pannonic woods with Quercus petraea and Carpinus betulus & yellow \\
	\hline 
	9110 & Luzulo-Fagetum beech forests
 & purple \\
	\hline
	\hline
\end{tabular}
    \caption{Natura 2000 forest habitats used in the classification.}
    \label{tab:habitats}
\end{center}
\end{table}

There were 125 areas segmented in red, green, blue and near-infrared channels of Sentinel-2 data \cite{ESASentinel}. We denote the segmented areas as $S_i$, where $i=1,\dots, N_S$, and $N_S=125$. This input from field experts contains 30 segmented areas of 91E0 habitat, which consists of mixed ash-alder alluvial forests in temperate and boreal Europe (Alno-padion, Alnion incanae, and Salicion albae); 29 91F0 areas, which consist of riparian mixed forests of Quercus robur, Ulmus laevis and Ulmus minor, Fraxinus excelsio or Fraxinus angustifolia, along the great rivers of the Atlantic and Middle European provinces; 32 91G0 areas, which consist of Pannonic woods with Quercus petrea and Carpinus betulus; and 34 9110 areas, which consist of Luzulo-Fagetum beech forests habitats. The Sentinel-2 data from September 10, 2018 covering Western Slovakia were used, as shown in the left square in Fig. \ref{fig:SK}. More specifically, we studied the Natura 2000 habitats in the Podunajsk\'{a} n\'{i}\v{z}ina lowland along the Danube river (see Fig. \ref{fig:habitats} left up), the Z\'{a}horská n\'{i}\v{z}ina lowland along the Morava river (see Fig. \ref{fig:habitats} right up) and the habitats in the Mal\'{e} Karpaty Mts. (see Fig. \ref{fig:habitats} bottom row).

\subsection{Obtaining of multispectral data and habitat spectral characteristics}

Sentinel-2 is a satellite of the European Space Agency (ESA) designed for the Copernicus European Union's Earth observation program focused on the observation of the atmosphere, land, seas and climate on Earth \cite{Copernicus}. All data acquired by the Sentinel-2 satellite are systematically processed, and only the products of this process are available for users \cite{ESASentinel}. We use the Level-2A product that provides Bottom Of Atmosphere (BOA) reflectance images and offers 17 channels. In addition to these 17 channels, we calculate one more, the normalized difference vegetation index (NDVI) \cite{EOS} giving further useful information on habitat status. Thus, for the classification and the feature space construction, we use 18 channels; and for every channel, we compute the statistical characteristics - the mean, the standard deviation, the minimum value and the maximum value in a prescribed image subarea $A$. Therefore, the feature space is the 72-dimensional Euclidean space; i.e., its dimensions are $k=72$.

\begin{figure}
 \begin{center}
	\includegraphics[width = 0.48\linewidth]{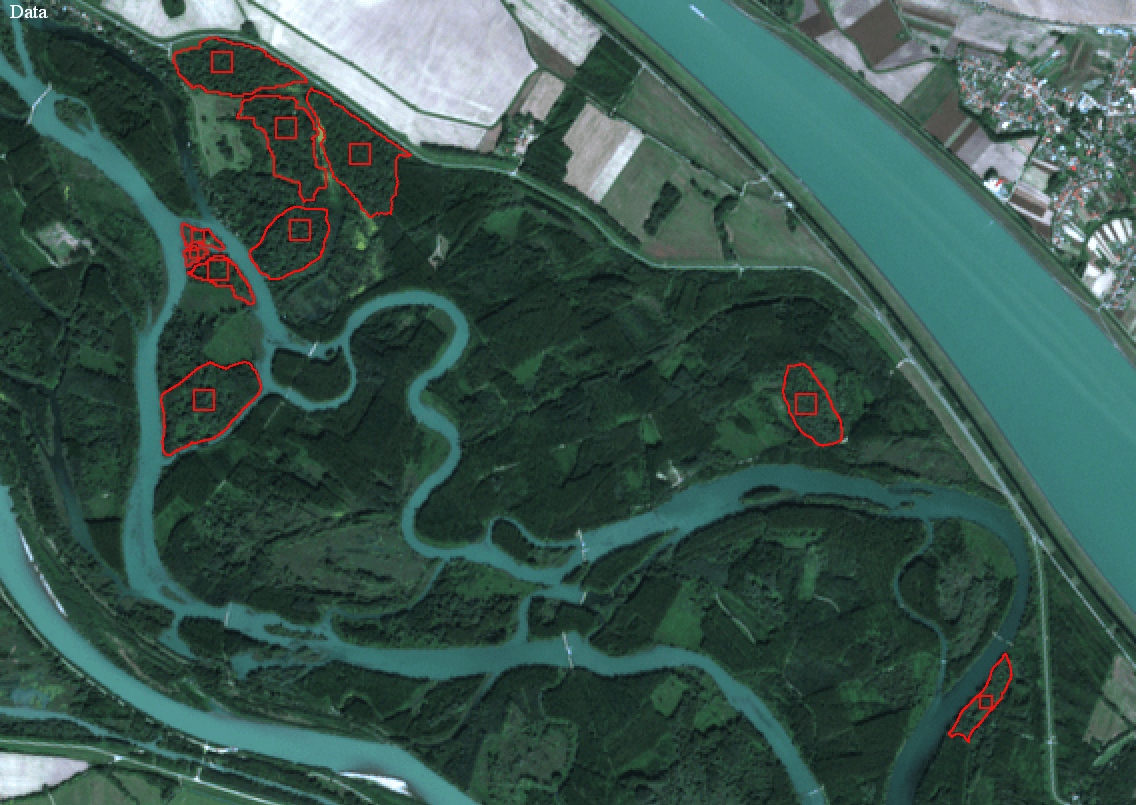}~~
	\includegraphics[width = 0.48\linewidth]{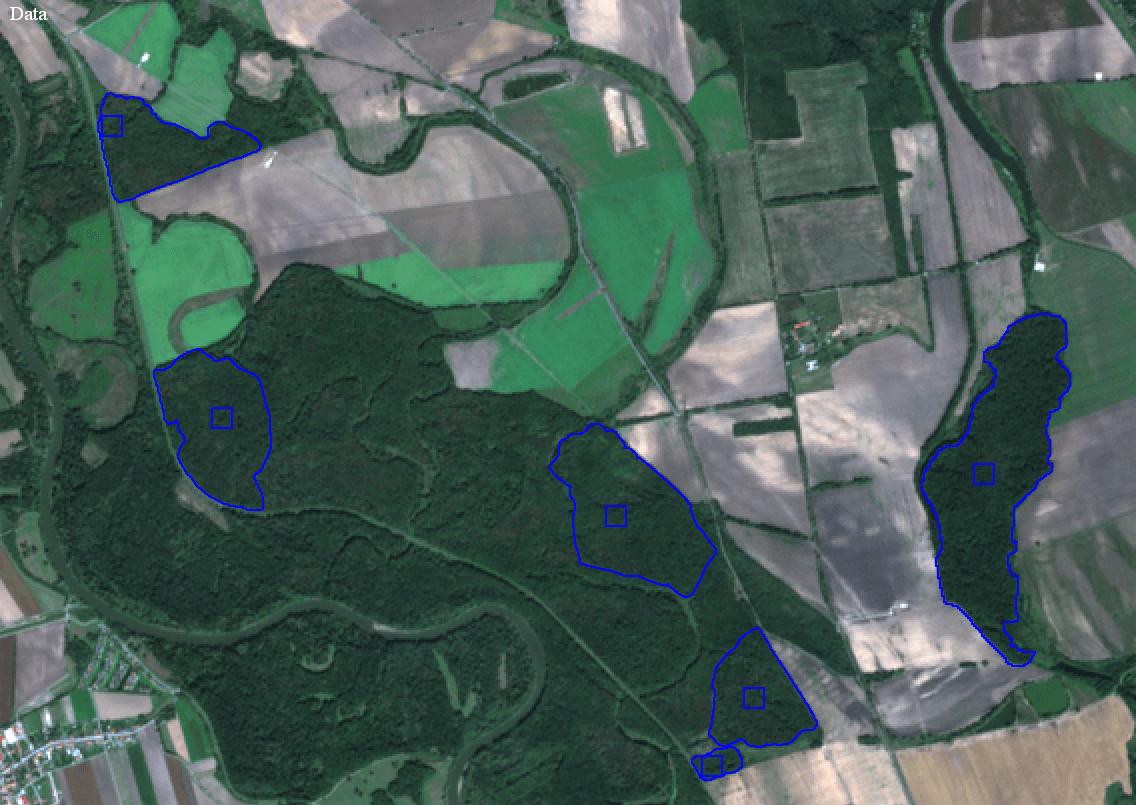} \\
	~\includegraphics[width = 0.48\linewidth]{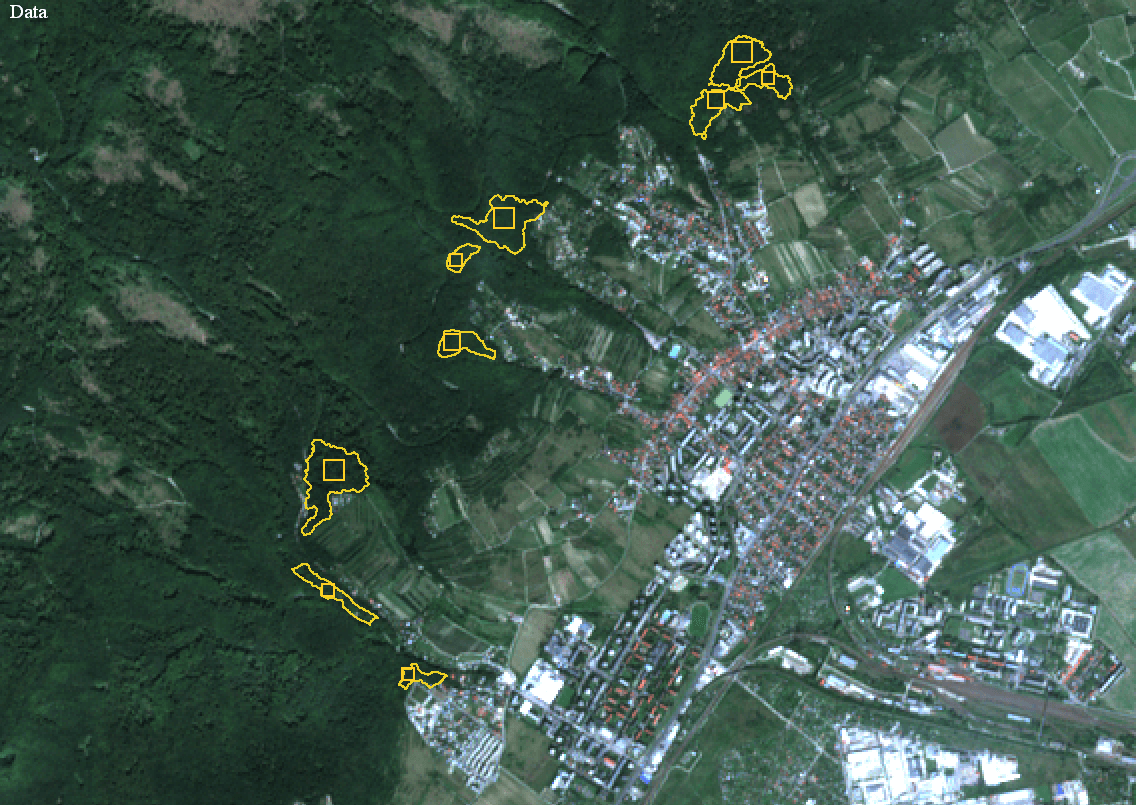}~~
	\includegraphics[width = 0.48\linewidth]{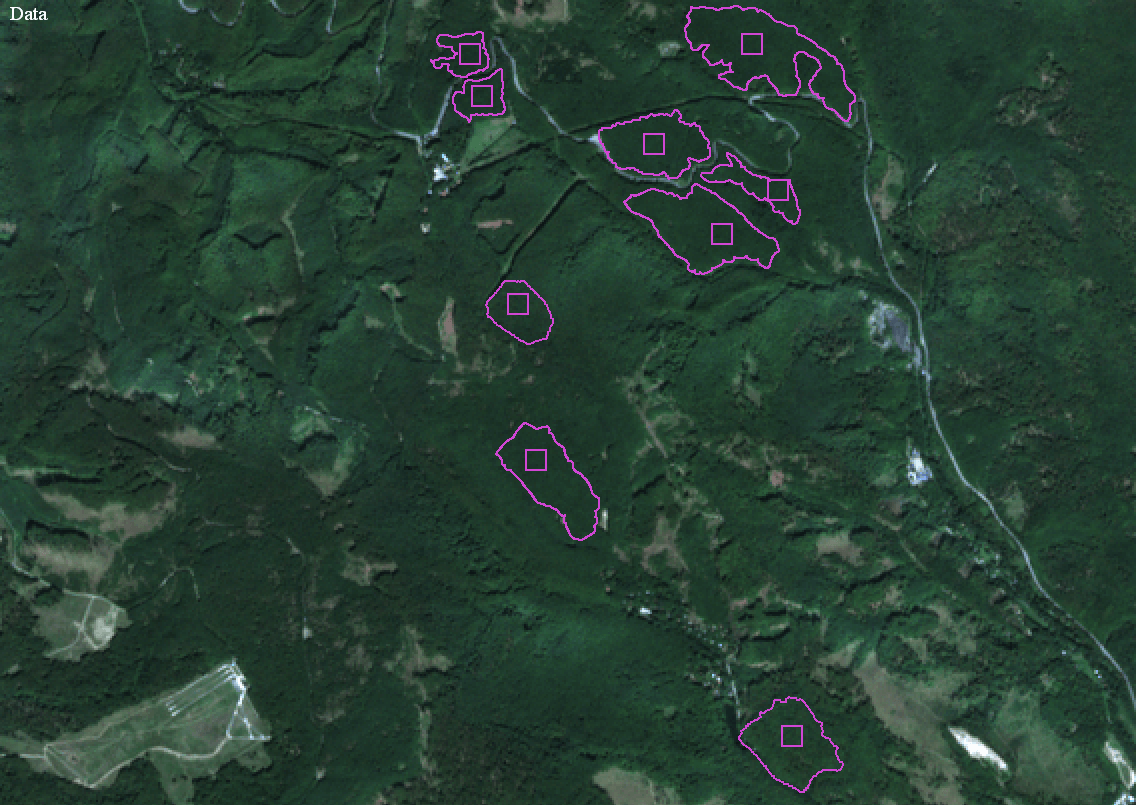}
	\caption{The subregions of Western Slovakia with segmented areas of protected Natura 2000 habitats. The segmented areas of {\rm 91E0} habitat are plotted by red curves, {\rm 91F0} habitat by blue curves, {\rm 91G0} habitat by yellow curves and {\rm 9110} habitat by purple curves. Together with the segmented areas we plot also squares inside the segmented areas which are used in the network learning phase and in the relevancy maps construction.}
	\label{fig:habitats_s}
 \end{center}
\end{figure}

For classification purposes, first, in each segmented area $S_i$, where $i=1,\dots, N_S$, we created a square $A_i=A(p_i,r)$ with a randomly chosen center in a pixel $p_i\in S_i$ and Chebyshev radius $r$. For large segmented areas, $r=5$ can be chosen; and for small areas, $r$ can be chosen smaller, as shown in Fig. \ref{fig:habitats_s}. Such squares are used for building the learning datasets and also for the construction of the so-called relevancy maps defined below. The statistical characteristics of the above mentioned 18 channels are computed for every square $A_i$, and they form the coordinates of points in the 72-dimensional feature space. The initial network graph $G$ is constructed such that every vertex of the graph $G$ is given by one such point corresponding to one square $A_i$ (or to one pixel $p_i$ or one segmented area $S_i$, as we may say).  

Since the feature space is high-dimensional, we reduce the dimensions. To detect and retain the maximal variance in the data, we apply Principal Component Analysis (PCA) \cite{wikiPCA, PCA}. We observed that in our application, the first two principal components are sufficient to describe the data variance; and thus, the dimension can be reduced from $k = 72$ to $k = 2$, which is simultaneously computationally tractable and yields convincing results.
	
\subsection{Relevancy map}
\label{subsec:relMap}
The relevancy map is a grayscale image with the same size as the images from Sentinel-2 optical channels. After finding the optimal parameters of the natural network, the square $A(p,r)$ is created in every image pixel $p$ (and not only inside the segmented areas as described above). For every $p$, the statistical characteristics of the square $A(p,r)$ are computed and considered as new observation $w(p)$. We note that to create required squares along the Sentinel-2 image boundary and compute the statistical characteristics there we use the reflection of corresponding values from the image interior to the exterior. The new observation $w(p)$ is added to the graph $ G $ as a new vertex. Every new observation $w(p)$ is classified by the natural network and its relevancy coefficient $R(w(p))$ is computed by (\ref{eq:relevancy1}). Finally, depending on the Chebyshev radius $r$ of the square $A(p,r)$, the relevancy map $ M_i^r $, $ i =1,\dots,N_C $, is defined for every cluster $C_i$ in every pixel $p$ as follows 
\begin{equation}
\begin{aligned}
	M_i^r(p) &= R(w(p)), \quad\ \ {\rm if} \ w(p) \ {\rm is\ classified\ into }\ C_i, \\
	M_i^r(p) &= 0, \quad\quad\quad\quad\quad\ {\rm if} \ w(p) \ {\rm is\ not\ classified\ into }\ C_i.
\end{aligned}
\label{eq:relevancymap}
\end{equation}

\subsection{Learning phase and network graph topology optimization}
First, let us consider the graph $G$ having $N_V=125$ vertices as described above and call it the (initial) learning dataset LDS125. 
All data from the learning dataset are labelled by the number of clusters to which they belong. As we stated above, we apply PCA to the original 72-dimensional points representing the data. PCA finds the coordinate system (basis) in which the highest variance (variability) of the data is in the first coordinate and it decreases subsequently with further coordinates. The change of the basis is represented by the linear transformation (matrix) that is applied to every point of the dataset, and we obtain coordinates of every point in the new coordinate system. Then, we are able to reduce the dimension of the feature space to $k=2$, considering only the first two coordinates of every point, as shown in the top left of Fig. \ref{fig:125x72_Squares_rev1}. As we observed experimentally, further coordinates do not help differentiate clusters and can be abandoned. After PCA matrix transformation, we also scale the coordinates of the points into the range 
$[0,1]\times [0,1]$ which helps in the model parameter tuning process. In fact, for any data, we can then use the same ranges of the natural network model parameters depending on the point distances.  
 
The main goals of the learning phase and the network graph optimization are to tune the parameters of the model (\ref{eq:kSysLE})-(\ref{eq:difKoefNum_newcom}) and optimize the structure of the graph $G$ itself to achieve the highest possible classification accuracy for observations from the learning dataset. To achieve that goal, we subsequently remove the cluster label from each vertex of graph $G$, set it as the new observation and classify it using the model, as shown in Fig. \ref{fig:125x72_Squares_rev1_n0}. We vary the model parameters $K_1$, $K_2$ and $\delta$ and choose the combination of parameters that results in the greatest number ($N_B$) of correctly classified observations from the learning dataset. If $N_V$ is the number of all observations, our goal is to achieve a success rate $N_B/N_V$ close to 1.

In the tuning process, the range for the diffusion coefficient parameters is $K_i \in [100, 5000]$ with the step size $K_{s_i} = 100 $, where $i = 1, 2$; and the range for parameter $\delta$ is the interval $[0.001, 0.1]$ with the step $\delta_s = 0.001$. As we tune three parameters only, we are able to go through all discrete parameters combinations for every new observation (a type of brute force approach) and find the combination that gives the highest possible classification accuracy $N_B$. For any parameter combination, we apply at most $n$ time steps of the natural network dynamics. We have chosen $n = 200$, but in most cases, the histogram stopping criterion is fulfilled for fewer time steps and classification is fast. For the numerical solution of the linear systems of equations (\ref{eq:kSysLE}), we use the SOR (Successive Overrelaxation) iterative method.
The first row of Table \ref{tab:LDS} shows the results. We achieved the best success rate $105/125$ by using the model parameters $K_1 = 2800$, $K_2 = 4700$ and $\delta = 0.004$.

\begin{table}[H]
\begin{center} 
\begin{tabular}{ c c c c c }
	\hline  
	Dataset name & \begin{tabular}{@{}c@{}} Correctly  \\ 
											classified
					\end{tabular} & 
					\begin{tabular}{@{}c@{}} Incorrectly  \\ 
											 classified
					\end{tabular} & 
											 Outliers &
											 Success rate \\
	\hline
	\hline 
	LDS125 & $ 105 $ & $ 16 $ & $ 4 $ & $ 84 \% $ \\
	LDS125adj & $ 113 $ & $ 10 $ & $ 2 $ & $ 90.4 \% $ \\
	LDS118 & $ 117 $ & $ 0 $ & $ 1 $ & $ 99.15 \% $ \\
	\hline
	\hline
\end{tabular}
    \caption{The results of the learning phase on datasets {\rm LDS125}, {\rm LDS125adj} and {\rm LDS118}.}
    \label{tab:LDS}
\end{center} 
\end{table}
	
The achieved the best success rate of 0.84 should be improved, so we performed further steps in the learning phase. It is quite clear that the initial random choice of the representative squares $A(p_i,r)$ inside the segmented areas $S_i$, where $i=1,\dots, N_S$, is not the optimal approach. The statistical characteristics of those squares do not necessarily correspond optimally to the statistical optical characteristics of the corresponding habitat in the Sentinel-2 image data. A strategy to solve this problem is based on spatial adjustment (shift) of representative squares inside the segmented areas $S_i$ such that we obtain a higher classification success rate with new (updated) vertices of the graph $G$. We can understand this step as a learning dataset adjustment. Since the representative squares have various radii (depending on the size of the segmented area), we construct the relevancy maps $M_i^r$, as shown in section \ref{subsec:relMap}, for $r=3, 4, 5$. Since the relevancy map gives the relevancy of classification for every image pixel, we will be able to compare the relevancies of the pixels inside the segmented areas and choose the pixel with the highest relevancy. The relevancy maps are constructed by using the graph $G$ corresponding to the initial learning dataset LDS125 and by using the optimal network parameters $K_1 = 2800$, $K_2 = 4700$, and $\delta = 0.004$. Then, we check whether there is an $r$ such that we are able to find a new pixel $p \in S_i$ in which $M_a^r(p)>>0$ (ideally close to 1), where $M_a^r$ is the relevancy map of the cluster to which the segmented area $S_i$ belongs. If we find such a pixel $p$ in the relevancy map $M_a^r$, the new square $A(p,r)$ is constructed. Every $A(p_i,r)$ for which it was possible to find the new representative square, with higher $M_a^r(p_i)$, is replaced by the new one and the adjusted learning dataset LDS125adj is created. The network parameter tuning process is conduct again, but it now uses the LDS125adj learning dataset. The result is shown in the second row of Tab. \ref{tab:LDS}; and the best success rate was $113/125=0.904$, which is higher than the previous one. It was achieved by using the parameters $K_1 = 4600$, $K_2 = 1700$, and $\delta = 0.002$. We can conclude that the adjustment of representative squares significantly increased the classification success rate.  

However, in the previous adjustment step, it was not always possible to find an adjusted representative square for all $S_i$, where $i=1,\dots, N_S$. This was caused by the fact that for some segmented areas $S_i$, only zero values of relevancy were computed in all interior pixels in all relevancy maps $M_a^r$, where $r=3, 4, 5$. We checked whether these only zero values in all relevancy maps were changed in those segmented areas when using LDS125adj dataset with its optimal parameters. Thus, we again constructed the relevancy maps. Checking the results, we found that there are still seven areas $S_i$ for which no inner pixel can be found such that it has relevancy greater than zero. It is clear that those areas cannot contribute in any way to increase the classification success rate. We call the area $S_i$ fulfilling the condition
\begin{equation}
	\label{eq:unclassifiable}
    M_a^r(p) = 0, \quad \forall p \in S_i, r=3, 4, 5,
\end{equation}	
unclassifiable and we remove all squares $A(p,r)$ corresponding to the unclassifiable areas from the learning dataset LDS125adj and create the final learning dataset LDS118 containing only $N_V=118$ labeled observations. Using this adjustment, we changed the topology of the graph $G$ itself, and we call this step the network graph topology optimization. The model parameters were tuned again and the best success rate of classification was $117/118 = 0.9915$, as shown in the third row of Tab. \ref{tab:LDS}, for the optimal parameters $K_1 = 3100$, $K_2 = 1500$, and $\delta = 0.003$. This success rate is high and allow us to use such optimally tuned (trained) natural network in practical applications presented in the next subsections.    

\section{Results and discussion}
\subsection{Relevancy maps for Western Slovakia}

\begin{figure}
 \begin{center}
	\includegraphics[width = 0.48\linewidth]{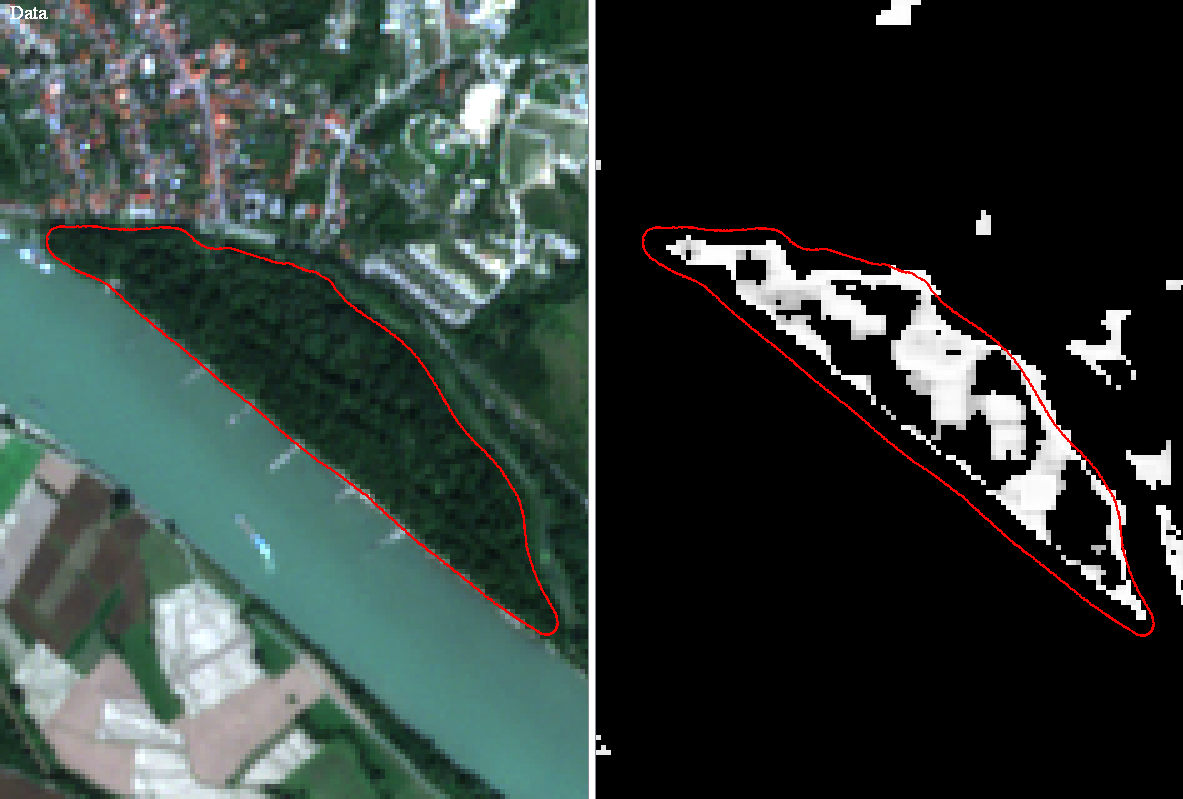}~~
	\includegraphics[width = 0.48\linewidth]{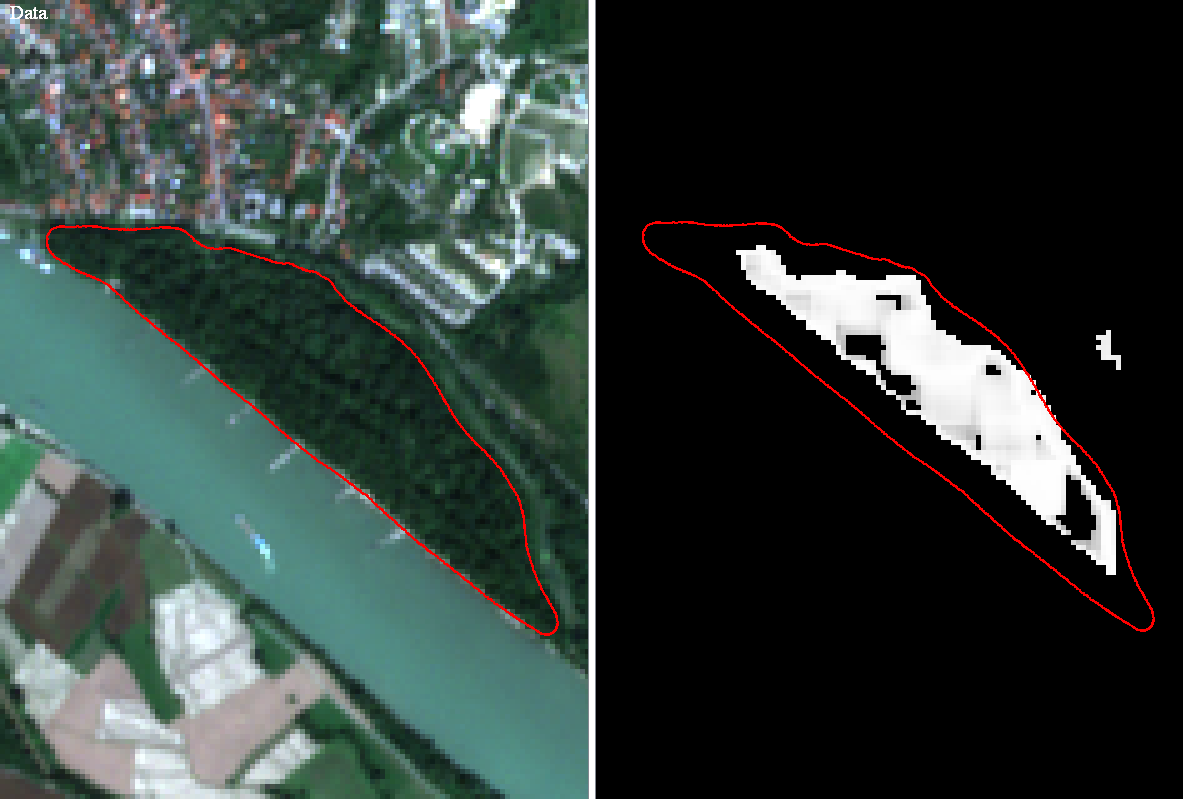} \\
	~\includegraphics[width = 0.48\linewidth]{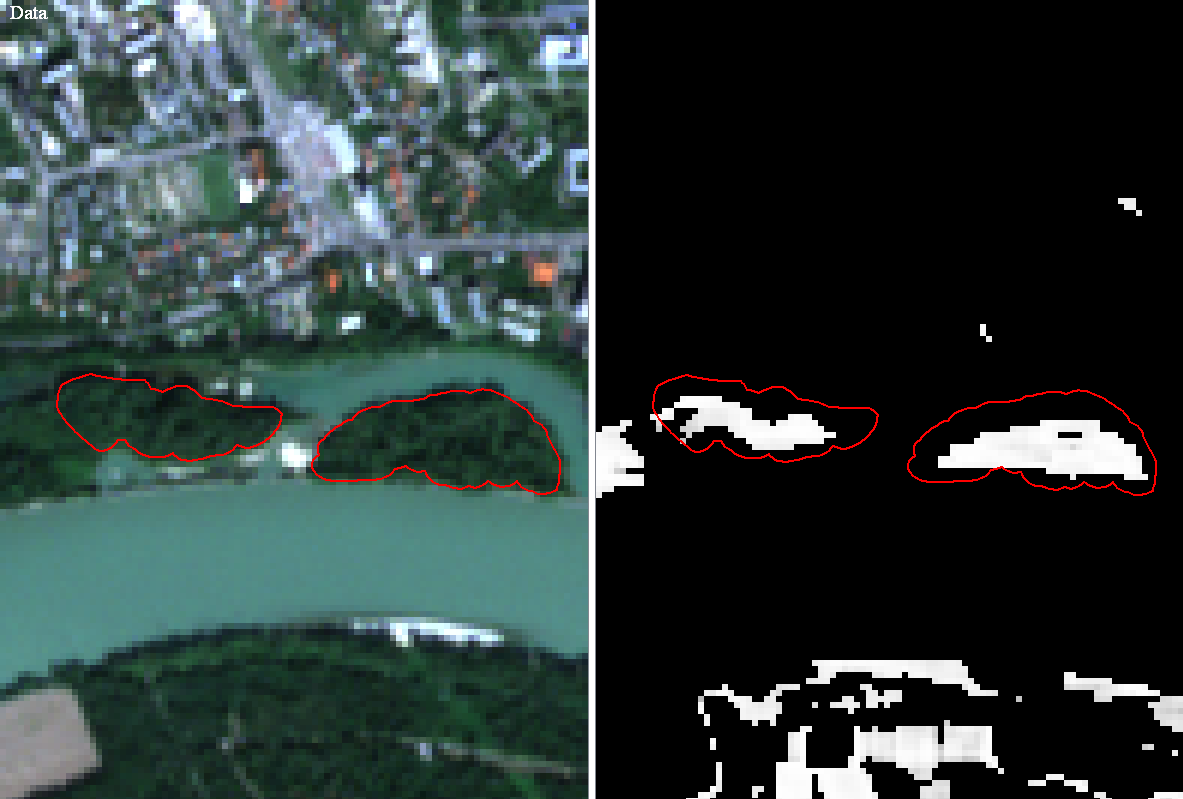}~~
	\includegraphics[width = 0.48\linewidth]{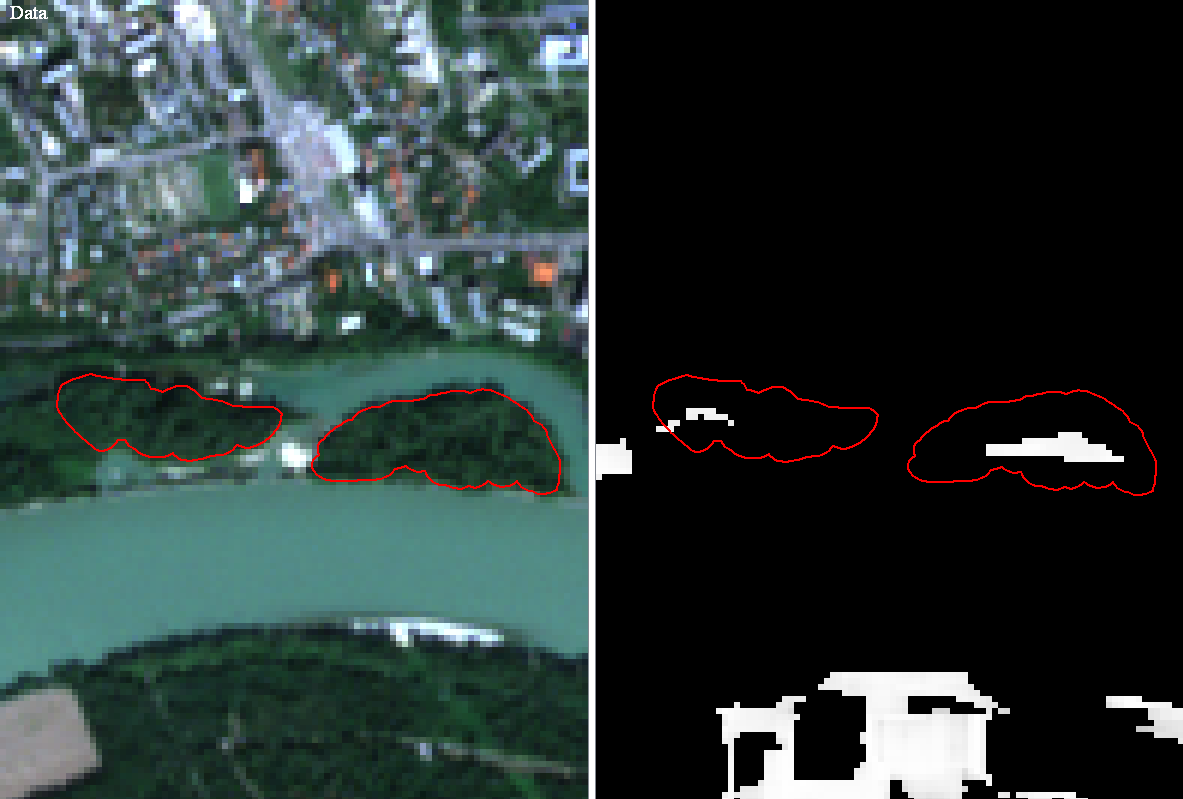}
	\caption{Each pair of pictures consists of the Sentinel-2 image (left part) and the relevancy map (right part). A large segmented area is depicted on the top pictures and two small segmented areas are depicted on the bottom pictures. Pictures on the left side depict the relevancy map obtained using the squares $A(p,3)$ and pictures on the right side depict the relevancy map with squares $A(p,5)$.}
	\label{fig:smallLarge}
 \end{center}
\end{figure}

In the Methods section, we discussed the learning process that led to the successfully optimized and trained natural network. Finally, the graph $G$ of the trained natural network contains $N_V=118$ vertices and the optimal parameters are $K_1 = 3100$, $K_2 = 1500$, and $\delta = 0.003$. This subsection is devoted to the construction of the final relevancy map for each explored Natura 2000 habitat. It was clear from the learning phase that the radius of the chosen square areas should vary between 3 (for small segmented areas) and 5 (for large segmented areas) to get the optimal results. For this reason, we compute three relevancy maps for the squares with radii $r=3$, $r=4$ and $r=5$, as shown in Fig. \ref{fig:smallLarge}. The final relevancy map $M_i^f$, where $i=1,\dots,N_C$, is obtained by taking the maximum of those three relevancies in every pixel $p$, i.e., we define it by
\begin{equation}
	M_i^f(p) = \max\limits_{r} M_i^r(p).
	\label{eq:finalrelevancy}
\end{equation}

\begin{figure}
 \begin{center}
	\includegraphics[width = 0.9\linewidth]{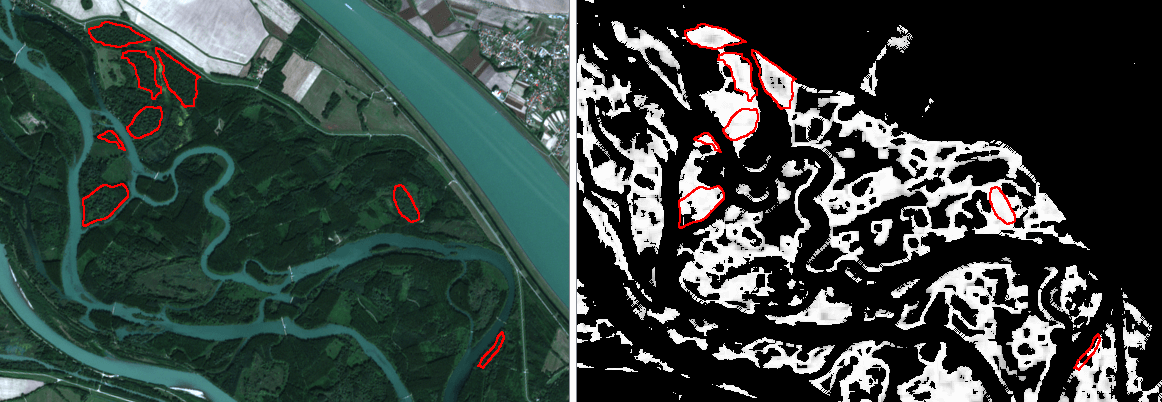} \\
	~\includegraphics[width = 0.9\linewidth]{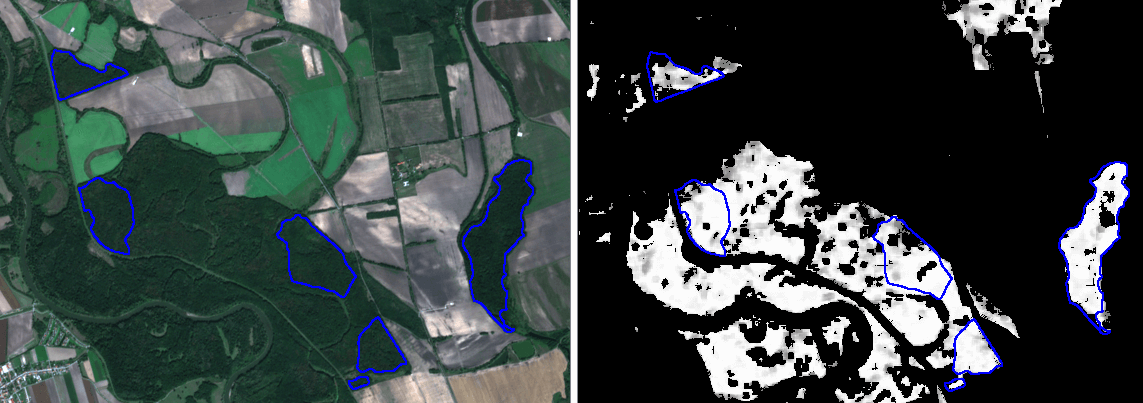}
	\caption{The segmented areas of {\rm 91E0} habitat (top) and {\rm 91F0} habitat (bottom) plotted on the Sentinel-2 image (left) and on the final relevancy map (right).}
	\label{fig:smallLarge_LD}
 \end{center}
\end{figure}

The qualitative (visual) comparisons of the Natura 2000 habitat segmented areas with the final relevancy maps for habitats 91E0 and 91F0 are plotted in Fig. \ref{fig:smallLarge_LD}. As the figure shows, the interior of the segmented areas on each relevancy map contains bright colors. This reflects the correct high relevancy of classification inside the segmented areas and thus the correct assignment of the image pixels to the Natura 2000 habitat. 
	
For the quantitative evaluation of the classification accuracy, we calculate the mean relevancy inside each segmented area. It is clear that the relevancy values inside the inner narrow band of width $r$ are irrelevant because the squares $A(p,r)$, for $p$ from the narrow band, also partially cover the pixels outside the area, cf. Fig. \ref{fig:smallLarge}. Thus, we shrink the boundary curve by the distance $r=3$, and the mean value of the relevancy is computed inside the smaller curve. In Tables \ref{tab:meanRCE0} - \ref{tab:meanRC10}, we show the mean relevancy inside the segmented areas of the 91E0, 91F0, 91G0 and 9110 Natura 2000 habitats in each of the final relevancy maps. 

The studied Natura 2000 habitat segmented areas obtain nonzero values in not only one final relevancy map. This occurs because the final relevancy map is constructed by the combination of the relevancy maps with different $r$s and by the fact that pixels inside the segmented areas can be classified differently by the natural network. We notice that most segmented areas obtain the highest mean relevancy for the habitat to which they belong. This is a feature supporting the usability of the final relevancy maps in the classification of new areas and finding new appearances of habitats. There are also some cases where the segmented area has the highest mean relevancy in a different cluster than that to which it should belong. This is, however, in a large majority of cases explainable. E.g., it may be caused by a similar species composition of habitats or by the very often difficult decision of experts in the field when classifying the habitat. Thus, we can conclude that the final relevancy maps give very useful information for an expert decision on habitat classification. This is further discussed below.  


The 91E0 softwood floodplain forests are the best detectable habitat among the studied habitats, as shown in Tab. \ref{tab:meanRCE0}, where all areas have the highest mean relevancy in the correct 91E0 cluster. The only disputable result is the habitat area 91E0\_Sihot, which was classified correctly as 91E0 but the relevancy values are high for both 91E0 and 91F0 habitats. It can be caused by the high cover of Ulmus minor, the hardwood tree species at this locality, which is also typical for the 91F0 habitat, but in combination with different trees. The floodplain forests at the Siho\v{t} Island are unique due to the specific co-occurrence of softwood floodplain forest species Populus sp. and Salix sp. together with Ulmus minor; therefore, the obtained relevancy values in Tab. \ref{tab:meanRCE0} are meaningful. It seems that when hardwood tree species are mixed in the 91E0 forest habitat, the mean relevancy of such an area also has a higher value for the 91F0 habitat cluster. This fact indicates recognition of transitional floodplain forest type that is common in natural conditions and also partly occurs in the localities in \v{C}\'{i}\v{c}ov (Cicov\_1, Cicov\_5) and Medve\v{d}ov (Medvedov\_2).

Within the 91F0 habitat, three areas are misclassified, and some others have close relevancies in both 91F0 and 91E0 habitats, as shown in Tab. \ref{tab:meanRCF0}. The segmented areas 91F0\_Kopac1, 91F0\_Vysoka\_pri\_Morave\_6 and 91F0\_Vysoka\_pri\_Morave\_8 were classified contrary to expert opinions. These forests represent a transitional type between 91E0 and 91F0 habitats, and experts classified them as 91F0. In many cases, it is hard to classify transitional floodplain forests fully objectively. The categorization is dependent on the subjective decision of the field expert. The relevancy shows possible classification within the 91E0 habitat of these three areas, which should be discussed among botany experts. In the alluvial country with natural characteristics, the softwood (91E0) and the hardwood (91F0) floodplain forests are often interconnected. They form transitional zones with hardwood tree species in the softwood habitat and, conversely, with softwood species in the hardwood floodplain forest. Our results show that the relevancy maps are sensitive to these transitional forests and can give information about such mixed composition of the floodplain forest and possible classification within the next-standing habitat.  

A few misclassified areas are found in the 91G0 habitat, as shown in Tab. \ref{tab:meanRCG0}. They are misclassified mainly within 91F0 instead of the correct 91G0. In addition, via a more in-depth look, we see that all misclassified areas contain Quercus petraea species together with Carpinus betulus. For the 91F0 habitat, Quercus robur is typical, but almost no difference between the two Quercus at the species level is identified. However, since the 91F0 habitat occurs in river alluvia while the 91G0 habitat occupies a lower hill, adding the information from a digital elevation model to the feature space may solve this problem.

In Table \ref{tab:meanRC10}, we present the results for the 9110 habitat areas. In some cases, they form transitional zones with 91G0 habitat because they share some tree species, Carpinus betulus, Tilia cordata, Prunus avium, Acer platanoides, etc. The beech forests in the Carpathians' upper parts are monodominant and thus are very easily detectable while forests in lower altitudes host more 91G0 species. Thus, the information from the digital elevation model can again improve the classification. One of the misclassified areas (9110\_Raca\_2) represents a very old-growth forest with a lower canopy cover, and probably some undergrowth species with different optical characteristics influenced the result.

\begin{table}[H]
\begin{center} \footnotesize
\begin{tabular}{ c c c c c }
	\hline  
	Area code & 91E0 habitat & 91F0 habitat & 91G0 habitat & 9110 habitat \\
	\hline
	\hline
	91E0\_Apalsky\_ostrov\_1    & 0.9163 & 0.0544 & 0 & 0 \\
    91E0\_Apalsky\_ostrov\_2    & 0.9149 & 0.2133 & 0 & 0 \\
    91E0\_Apalsky\_ostrov\_3    & 0.8549 & 0.2321 & 0 & 0 \\
    91E0\_Bodiky\_1             & 0.6991 & 0.0110 & 0 & 0 \\
    91E0\_Bodiky\_2             & 0.7892 & 0.1230 & 0 & 0 \\
    91E0\_Bodiky\_3             & 0.8103 & 0.2247 & 0 & 0 \\
    91E0\_Bodiky\_4             & 0.8419 & 0.1381 & 0 & 0 \\
    91E0\_Bodiky\_5             & 0.8567 & 0.0384 & 0 & 0 \\
    91E0\_Bodiky\_6             & 0.9435 & 0.0049 & 0 & 0 \\
    91E0\_Bodiky\_7             & 0.9632 & 0.1884 & 0 & 0 \\
    91E0\_Bodiky\_8             & 0.8816 & 0.1825 & 0 & 0 \\
    91E0\_Bodiky\_9             & 0.4281 & 0      & 0 & 0 \\
    91E0\_Cicov\_1              & 0.8141 & 0.4700 & 0 & 0 \\
    91E0\_Cicov\_2              & 0.9308 & 0.0598 & 0 & 0 \\
    91E0\_Cicov\_3              & 0.2660 & 0      & 0 & 0 \\
    91E0\_Cicov\_4              & 0.9637 & 0.3541 & 0 & 0 \\
    91E0\_Cicov\_5              & 0.8504 & 0.5472 & 0 & 0 \\
    91E0\_Karloveske\_rameno\_1 & 0.6677 & 0      & 0 & 0 \\
    91E0\_Karloveske\_rameno\_2 & 0.6313 & 0.0267 & 0 & 0 \\
    91E0\_Kollarovo             & 0.8372 & 0      & 0 & 0 \\
    91E0\_Medvedov\_1           & 0.7647 & 0.3731 & 0 & 0 \\
    91E0\_Medvedov\_2           & 0.8136 & 0.4445 & 0 & 0 \\
    91E0\_Morusovy\_luh         & 0.4209 & 0      & 0 & 0 \\
    91E0\_Sala                  & 0.5539 & 0      & 0 & 0 \\
    91E0\_Sihot                 & 0.5781 & 0.4892 & 0 & 0 \\
    91E0\_Slovansky\_ostrov     & 0.7400 & 0.3391 & 0 & 0 \\
    91E0\_Vahovce               & 0.4519 & 0.0124 & 0 & 0 \\
    91E0\_Velky\_lel            & 0.5693 & 0      & 0 & 0 \\
	\hline
	\hline
\end{tabular}
    \caption{The segmented areas of the Natura 2000 habitat {\rm 91E0} with the mean value of the relevancy in each final relevancy map.}
    \label{tab:meanRCE0}
\end{center} 
\end{table}

\begin{table}[H]
\begin{center} \footnotesize
\begin{tabular}{ c c c c c }
	\hline  
	Area code & 91E0 habitat & 91F0 habitat & 91G0 habitat & 9110 habitat \\
	\hline
	\hline
	91F0\_Bazantnica             & 0.0788 & 0.6781 & 0.3776 & 0.0013 \\
    91F0\_Bogdalicky\_vrch\_1    & 0.1757 & 0.8159 & 0.1694 & 0      \\
    91F0\_Bogdalicky\_vrch\_2    & 0.2456 & 0.7487 & 0.2778 & 0      \\
    91F0\_Bogdalicky\_vrch\_3    & 0.1584 & 0.8833 & 0.2594 & 0      \\
    91F0\_Bogdalicky\_vrch\_4    & 0.2511 & 0.9068 & 0.0460 & 0      \\
    91F0\_Brestovany\_1          & 0.0537 & 0.8097 & 0.4241 & 0      \\
    91F0\_Brestovany\_2          & 0.3178 & 0.6955 & 0.0457 & 0      \\
    91F0\_Brestovany\_3          & 0.2889 & 0.8301 & 0.2675 & 0      \\
    91F0\_Brestovany\_4          & 0.5076 & 0.8102 & 0.0035 & 0      \\
    91F0\_Dolnyles               & 0.1420 & 0.7016 & 0.4337 & 0      \\
    91F0\_Feldsky\_les\_1        & 0.2880 & 0.7741 & 0.0136 & 0      \\
    91F0\_Feldsky\_les\_2        & 0.1991 & 0.8693 & 0.0161 & 0      \\
    91F0\_Kopac\_1               & 0.8349 & 0.3646 & 0      & 0      \\
    91F0\_Kopac\_2               & 0.4908 & 0.6242 & 0.0421 & 0      \\
    91F0\_Kopac\_3               & 0.0526 & 0.9502 & 0      & 0      \\
    91F0\_Kopac\_4               & 0.3927 & 0.8416 & 0      & 0      \\
    91F0\_Suchohrad\_1           & 0.4707 & 0.7697 & 0.0675 & 0      \\
    91F0\_Suchohrad\_2           & 0.2324 & 0.6818 & 0.1897 & 0      \\
    91F0\_Suchohrad\_3           & 0.4620 & 0.5770 & 0.0400 & 0      \\
    91F0\_Vysoka\_pri\_Morave\_1 & 0.0858 & 0.7792 & 0.3458 & 0      \\
    91F0\_Vysoka\_pri\_Morave\_2 & 0.0693 & 0.7032 & 0.2892 & 0.0004 \\
    91F0\_Vysoka\_pri\_Morave\_3 & 0.1795 & 0.8809 & 0.0513 & 0      \\
    91F0\_Vysoka\_pri\_Morave\_4 & 0.1535 & 0.7929 & 0.3317 & 0      \\
    91F0\_Vysoka\_pri\_Morave\_5 & 0.3218 & 0.7093 & 0.4345 & 0      \\
    91F0\_Vysoka\_pri\_Morave\_6 & 0.7894 & 0.4994 & 0      & 0      \\
    91F0\_Vysoka\_pri\_Morave\_7 & 0.2338 & 0.7984 & 0.0859 & 0      \\
    91F0\_Vysoka\_pri\_Morave\_8 & 0.5009 & 0.4573 & 0.0640 & 0      \\
    91F0\_Vysoka\_pri\_Morave\_9 & 0.1236 & 0.7318 & 0.2276 & 0      \\    	\hline
	\hline
\end{tabular}
    \caption{The segmented areas of the Natura 2000 habitat {\rm 91F0} with the mean value of the relevancy in each final relevancy map.}
    \label{tab:meanRCF0}
\end{center} 
\end{table}

\begin{table}[H]
\begin{center} \footnotesize
\begin{tabular}{ c c c c c }
	\hline  
	Area code & 91E0 habitat & 91F0 habitat & 91G0 habitat & 9110 habitat \\
	\hline
	\hline
	91G0\_Casta\_1     & 0      & 0.0098 & 0.8076 & 0.2909 \\
    91G0\_Casta\_2     & 0      & 0.1127 & 0.7744 & 0.2064 \\
    91G0\_Dubova\_1    & 0.1572 & 0.1403 & 0.5290 & 0      \\
    91G0\_Dubova\_2    & 0.0334 & 0.1036 & 0.8082 & 0.0063 \\
    91G0\_Limbach\_1   & 0      & 0.2438 & 0.9026 & 0      \\
    91G0\_Limbach\_2   & 0.0385 & 0.7094 & 0.6549 & 0      \\
    91G0\_Limbach\_3   & 0      & 0      & 0.9935 & 0      \\
    91G0\_Limbach\_4   & 0      & 0.1618 & 0.9733 & 0      \\
    91G0\_Limbach\_5   & 0      & 0.0433 & 0.9933 & 0      \\
    91G0\_Losonec\_1   & 0.0913 & 0.3749 & 0.6612 & 0.0239 \\
    91G0\_Losonec\_2   & 0.0282 & 0.1941 & 0.8678 & 0.0339 \\
    91G0\_Losonec\_3   & 0      & 0.0109 & 0.8967 & 0.0587 \\
    91G0\_Pezinok\_1   & 0.0412 & 0.0428 & 0.1943 & 0.7222 \\
    91G0\_Pezinok\_2   & 0      & 0      & 0.9926 & 0      \\
    91G0\_Pezinok\_3   & 0      & 0.0089 & 0.9499 & 0.0809 \\
    91G0\_Pezinok\_4   & 0.0027 & 0.0702 & 0.9153 & 0.0633 \\
    91G0\_Pezinok\_5   & 0      & 0.2264 & 0.9227 & 0      \\
    91G0\_Pezinok\_6   & 0.0796 & 0.0354 & 0.7120 & 0.1048 \\
    91G0\_Raca\_1      & 0.1290 & 0.9346 & 0.0971 & 0      \\
    91G0\_Raca\_2      & 0      & 0      & 0.3056 & 0.8097 \\
    91G0\_Raca\_3      & 0.0734 & 0.2068 & 0.8630 & 0.0165 \\
    91G0\_Raca\_4      & 0      & 0.0260 & 0.8926 & 0      \\
    91G0\_Raca\_5      & 0.2952 & 0.3849 & 0.1116 & 0      \\
    91G0\_Raca\_6      & 0.0258 & 0.1099 & 0.9134 & 0.0050 \\
    91G0\_Raca\_7      & 0.3668 & 0.5320 & 0.1257 & 0      \\
    91G0\_Raca\_8      & 0      & 0.7130 & 0.7517 & 0      \\
    91G0\_Raca\_9      & 0.0085 & 0      & 0.9766 & 0      \\
    91G0\_Smolenice\_1 & 0.0164 & 0.0396 & 0.8445 & 0.1793 \\
    91G0\_Smolenice\_2 & 0      & 0.0000 & 0.9860 & 0      \\
    91G0\_Smolenice\_3 & 0      & 0.2037 & 0.8959 & 0      \\
    91G0\_Smolenice\_4 & 0      & 0      & 0.9758 & 0.0641 \\
	\hline
	\hline
\end{tabular}
    \caption{The segmented areas of the Natura 2000 habitat {\rm 91G0} with the mean value of the relevancy in each final relevancy map.}
    \label{tab:meanRCG0}
\end{center} 
\end{table}

\begin{table}[H]
\begin{center} \footnotesize
\begin{tabular}{ c c c c c }
	\hline  
	Area code & 91E0 habitat & 91F0 habitat & 91G0 habitat & 9110 habitat \\
	\hline
	\hline
	9110\_Borinka\_1        & 0.0417 & 0.0857 & 0.2324 & 0.5473 \\
    9110\_Borinka\_2        & 0      & 0      & 0.3016 & 0.7477 \\
    9110\_Borinka\_3        & 0.0041 & 0.0272 & 0.6455 & 0.4040 \\
    9110\_Limbach\_1        & 0      & 0.1060 & 0.8554 & 0.2335 \\
    9110\_Limbach\_2        & 0      & 0      & 0.0332 & 0.9803 \\
    9110\_Limbach\_3        & 0      & 0      & 0      & 0.5491 \\  
    9110\_Limbach\_4        & 0      & 0      & 0.0003 & 0.4209 \\  
    9110\_Limbach\_5        & 0      & 0      & 0.0848 & 0.8750 \\  
    9110\_Limbach\_6        & 0      & 0      & 0.1501 & 0.8651 \\  
    9110\_Limbach\_7        & 0.0087 & 0.0240 & 0.1488 & 0.8900 \\  
    9110\_Limbach\_8        & 0.0129 & 0.0378 & 0.2263 & 0.6683 \\  
    9110\_Limbach\_9        & 0.0253 & 0.0537 & 0.2309 & 0.6883 \\
    9110\_Limbach\_10       & 0      & 0.0000 & 0.0076 & 0.9610 \\
    9110\_Limbach\_11       & 0.0030 & 0.0035 & 0.0141 & 0.8516 \\
    9110\_Modra\_Piesok\_1  & 0.0465 & 0.0312 & 0.2866 & 0.7078 \\
    9110\_Modra\_Piesok\_2  & 0.0462 & 0.1169 & 0.2982 & 0.6277 \\
    9110\_Modra\_Piesok\_3  & 0      & 0.0000 & 0.0175 & 0.9893 \\
    9110\_Modra\_Piesok\_4  & 0.0413 & 0.0368 & 0.0485 & 0.9425 \\
    9110\_Pezinok\_1        & 0.0017 & 0.0007 & 0.0044 & 0.7455 \\
    9110\_Pezinok\_2        & 0      & 0.0087 & 0.0377 & 0.8132 \\
    9110\_Pezinska\_Baba\_1 & 0      & 0      & 0.0145 & 0.7932 \\
    9110\_Pezinska\_Baba\_2 & 0.0028 & 0.0011 & 0.0096 & 0.9372 \\
    9110\_Pezinska\_Baba\_3 & 0.0071 & 0.0127 & 0.2394 & 0.8769 \\
    9110\_Pezinska\_Baba\_4 & 0      & 0.0141 & 0.0661 & 0.9621 \\
    9110\_Pezinska\_Baba\_5 & 0.0093 & 0.0039 & 0.0842 & 0.9497 \\
    9110\_Raca\_1           & 0      & 0      & 0.0133 & 0.9364 \\
    9110\_Raca\_2           & 0.0034 & 0.2190 & 0.7639 & 0.0372 \\
    9110\_Raca\_3           & 0      & 0      & 0.0682 & 0.9609 \\
    9110\_Raca\_4           & 0.0402 & 0.1584 & 0.6640 & 0.3444 \\
    9110\_Raca\_5           & 0.0250 & 0.1000 & 0.2317 & 0.4714 \\
    9110\_Smolenice         & 0      & 0      & 0.1020 & 0.8774 \\
	\hline
	\hline
\end{tabular}
    \caption{The segmented areas of the Natura 2000 habitat {\rm 9110} with the mean value of the relevancy in each final relevancy map.}
    \label{tab:meanRC10}
\end{center} 
\end{table}

The natural vegetation has a continuous character, and Natura 2000 habitats stated in the habitats directive describe the essential, most characteristic forms of habitats. In real nature, these types are often interconnected and hard to classify, even by experts in the field. The presented results show the high classification success; moreover, it gives us information about the admixture of some tree species typical for next-standing habitats. This finding opens new opportunities in phytosociological and ecological research, especially in the fields of ecotone zones and vegetation gradients.

\subsection{Exploration of other regions in Slovakia and finding new appearances of protected habitats}
After obtaining the successfully trained natural network, we can explore further areas of the occurrence of protected habitats in territories neighboring Western Slovakia, as shown in Fig. \ref{fig:SK}. We explored the areas of Central and Eastern Slovakia and obtained very promising results. By using the relevancy maps, we confirmed the appearance of the 91F0 habitat around the Latorica River in the south of Eastern Slovakia and found the new appearance of 91E0 protected habitat along the Rimava River near the village of Dubovec in the southern part of Central Slovakia. The main goal of the relevancy map construction, to find the new appearances of protected habitats in an automated way, was thus achieved.
	
\begin{figure}
 \begin{center}
	\includegraphics[width = 0.48\linewidth]{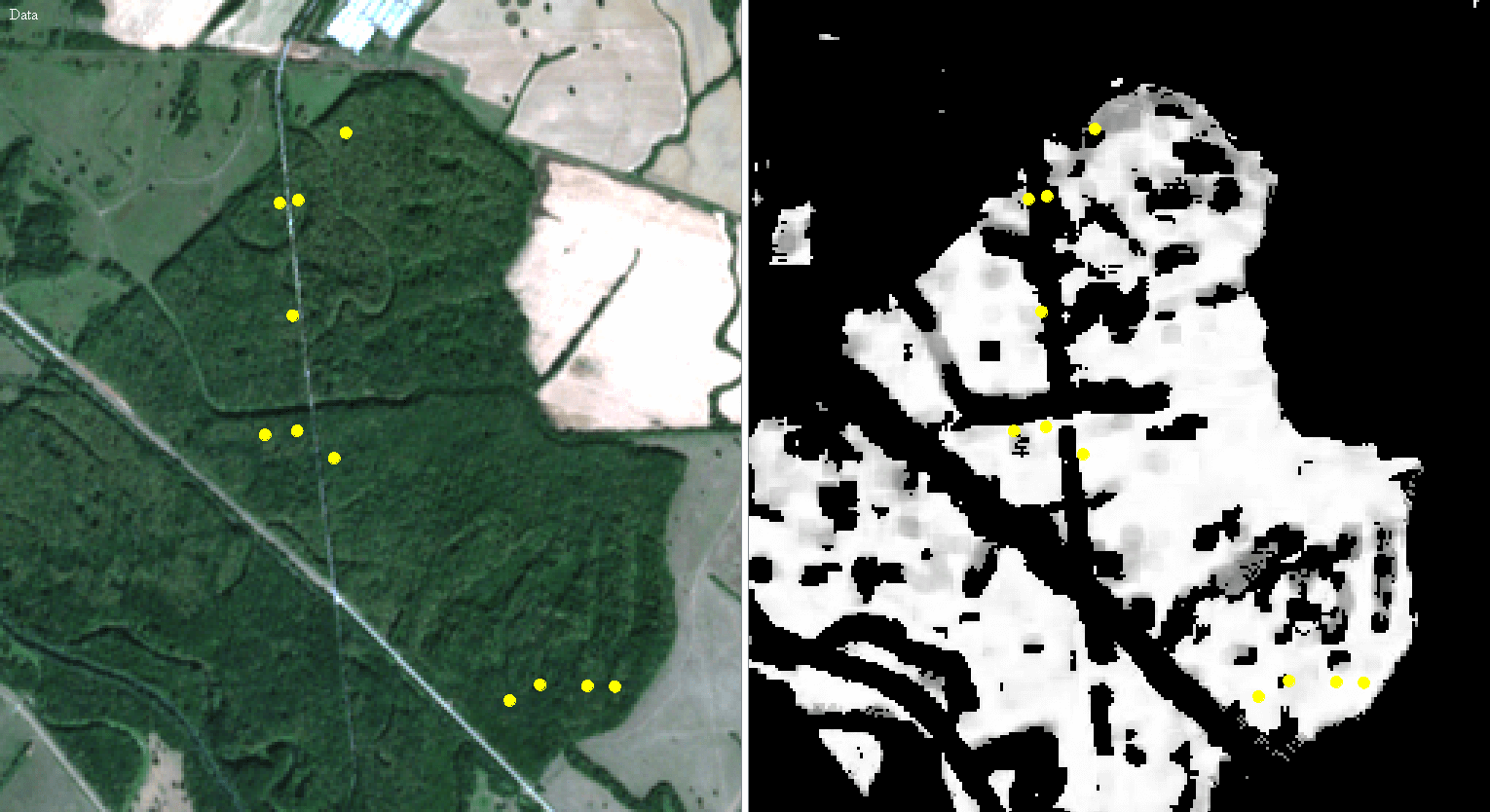}~~
	\includegraphics[width = 0.48\linewidth]{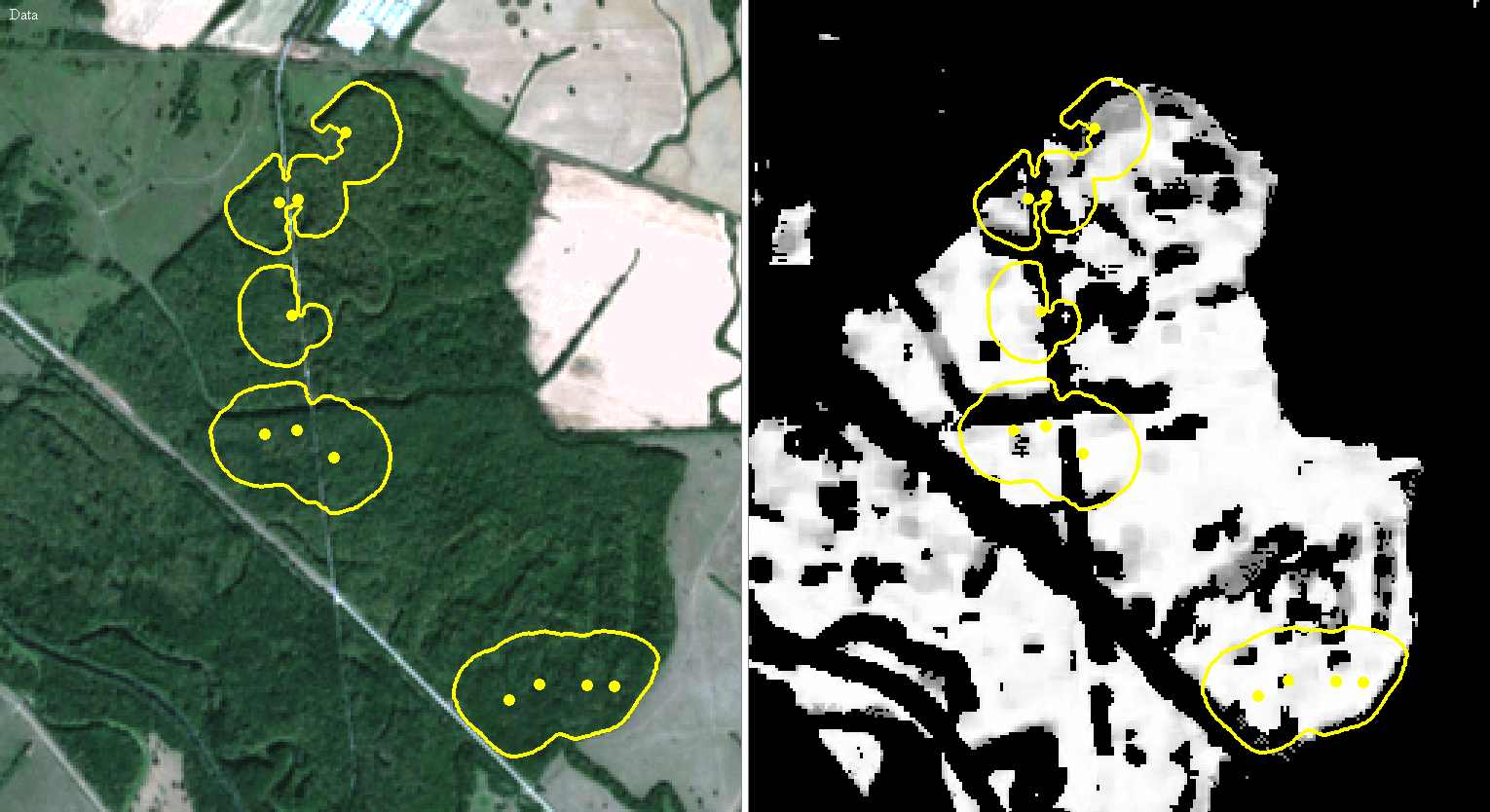} \\
	~\includegraphics[width = 0.48\linewidth]{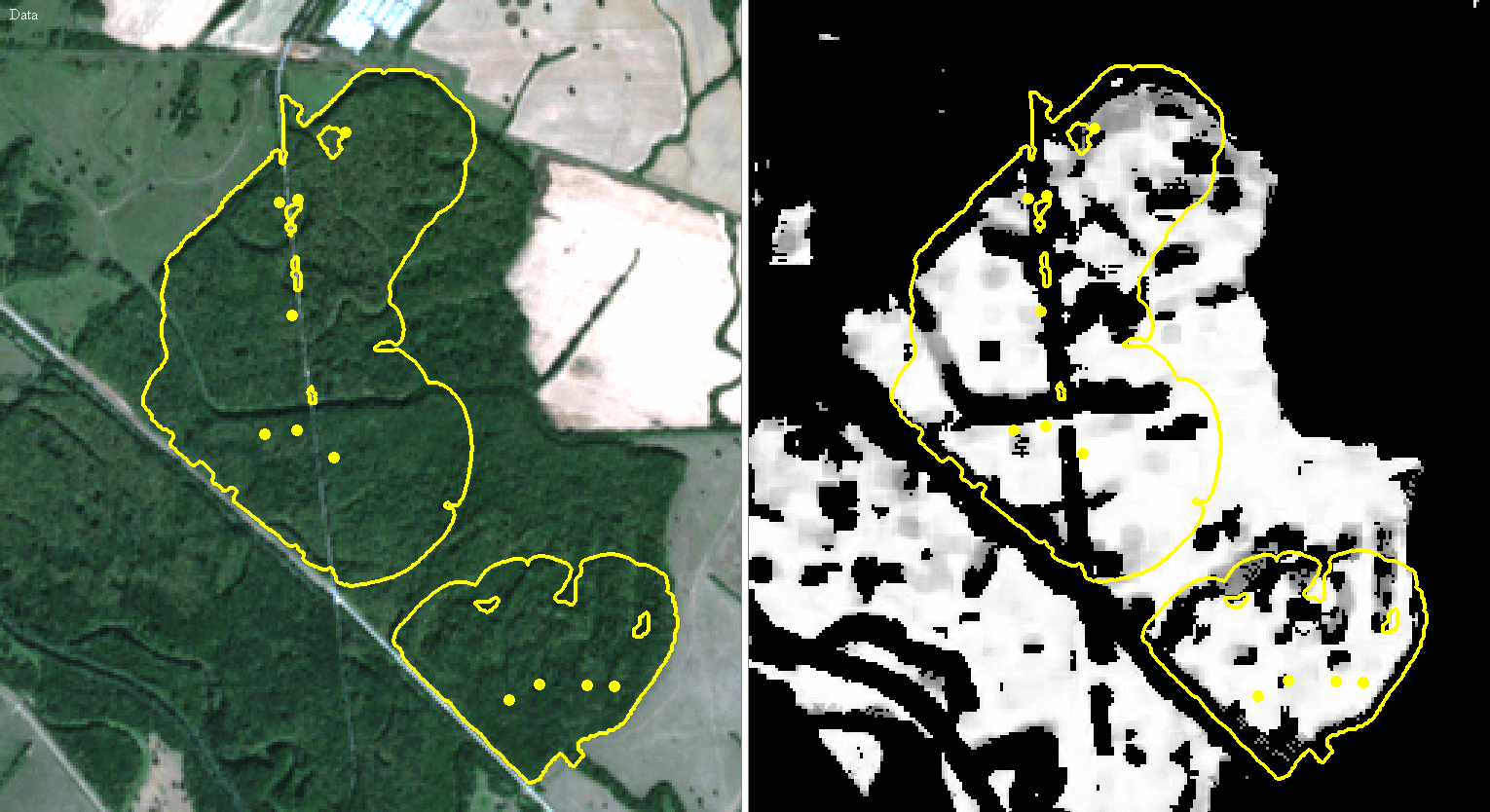}~~
	\includegraphics[width = 0.48\linewidth]{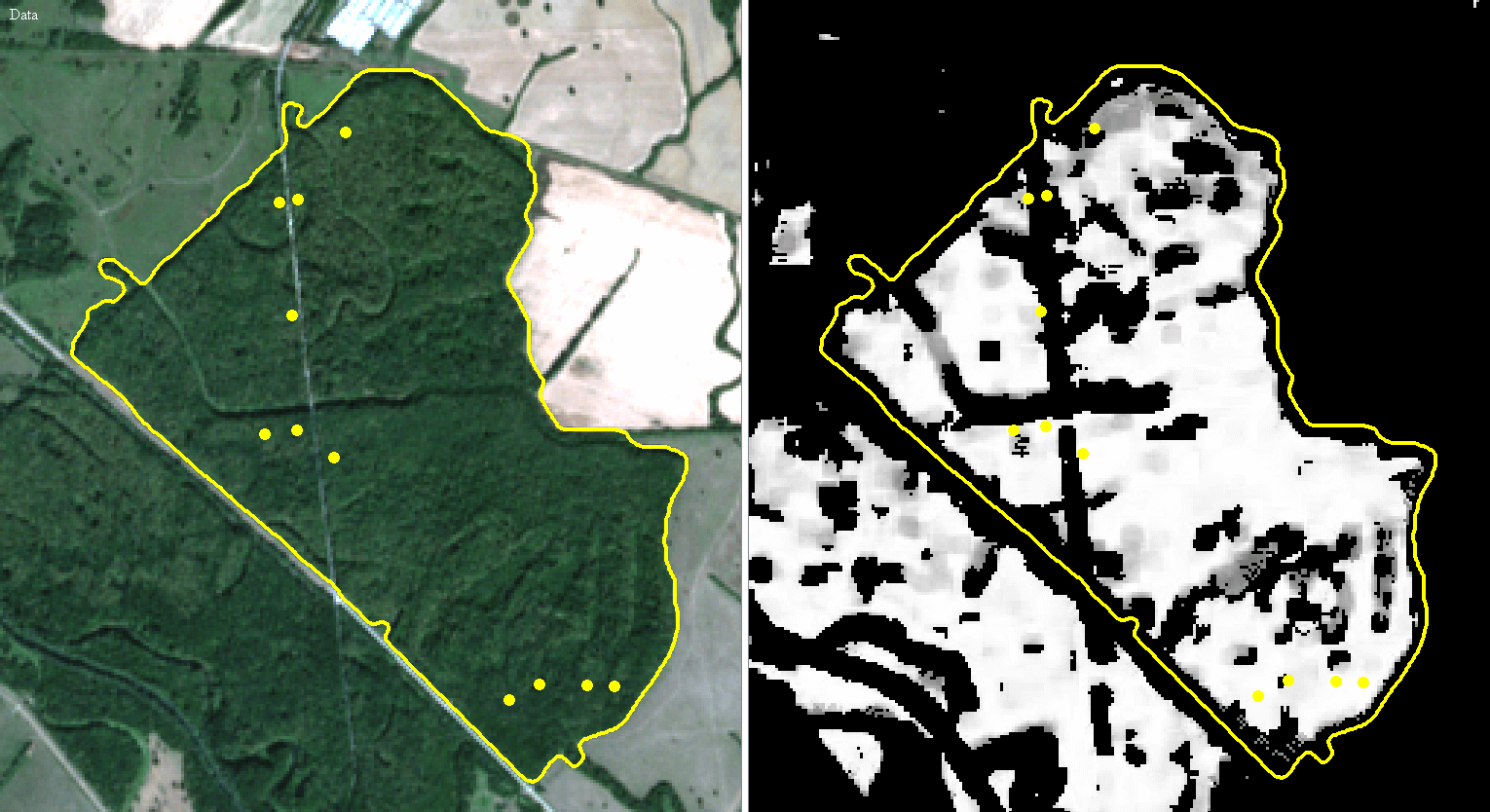}
	\caption{The area around the Latorica river in Eastern Slovakia, the Sentinel-2 image (left) and relevancy map (right) with subsequent automatic segmentation of {\rm 91F0} habitat area. }
	\label{fig:Latorica}
 \end{center}
\end{figure}

During the Natura 2000 mapping campaigns, the 91F0 floodplain forests around the Latorica River were sampled using vegetation plots. We use the trained natural network to create the final relevancy maps for the Eastern Slovakia region, and the details of the Sentinel-2 image and the 91F0 final relevancy map can be seen in Fig. \ref{fig:Latorica}. The figure shows 11 points (yellow color) representing the vegetation plots - phytosociological relevés, which mark the appearance of the habitat. Usually, they are not far from the boundary of the habitat. First, we see the correspondence of these points and bright colors in the relevancy map of the 91F0 habitat, which indicates the correct classification of the pixels along the Latorica River. Furthermore, we automatically segmented the habitat area \cite{segAuto} starting the segmentation from 11 habitat marking points. The evolution of the segmentation curves undergoing topological changes is presented in Fig. \ref{fig:Latorica} with the final segmentation result plotted in the bottom right picture both on the Sentinel-2 image and relevancy map. The automatic segmentation identified the compact habitat area and  also included a thin river branch, a forest road and small areas with younger floodplain forests, which is accepted because these objects are common parts of the habitat. Furthermore, the relevancy map was able to detect even these small parts, and high relevancy values were assigned only to the best representative area of the habitat. We computed also the mean value of the relevancy in the 91F0 relevancy map inside the segmented area; and it is $ 0.6548 $, which indicates the correct classification of the area along the Latorica River.

\begin{figure}
 \begin{center}
	\includegraphics[width = 0.48\linewidth]{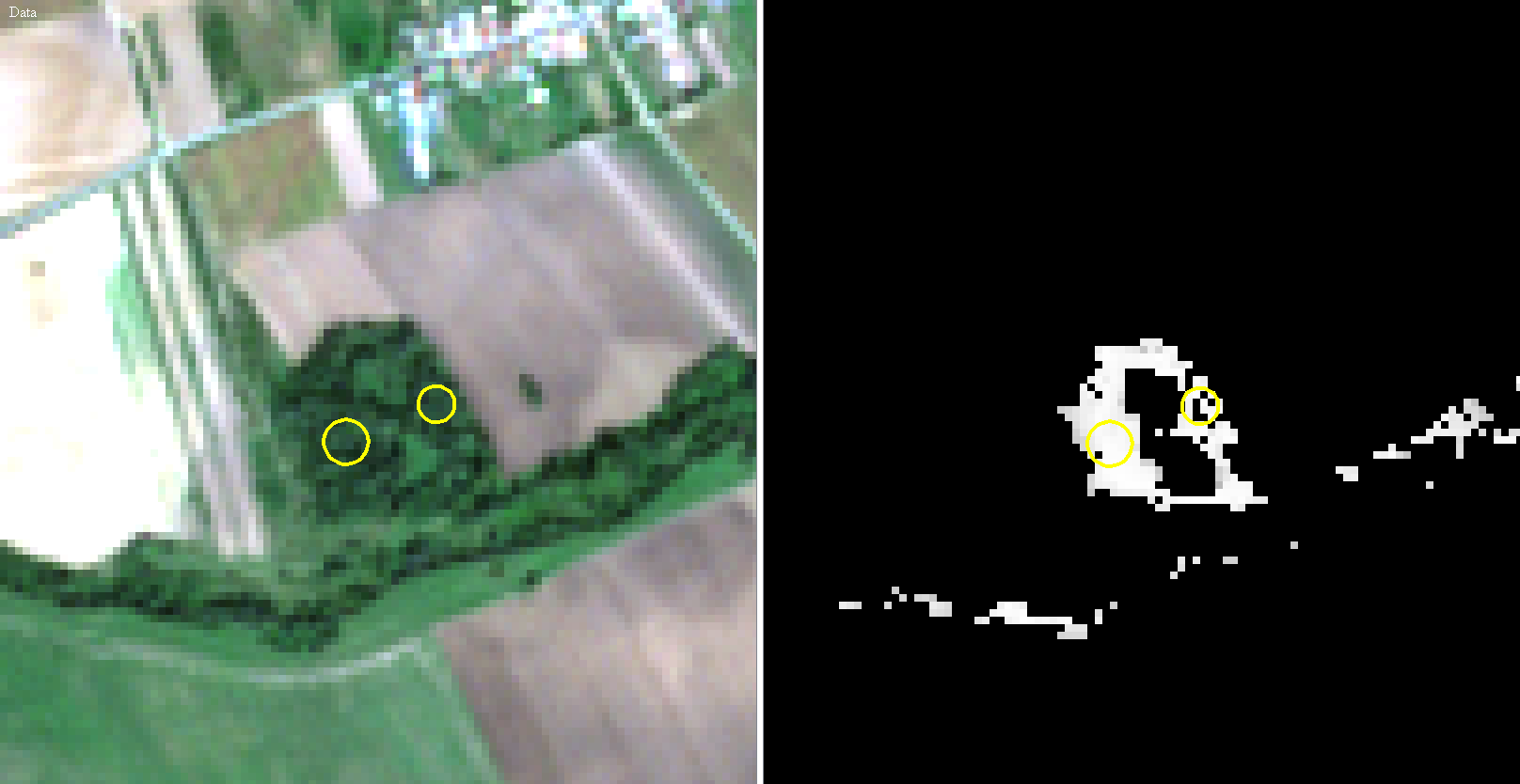}~~
	\includegraphics[width = 0.48\linewidth]{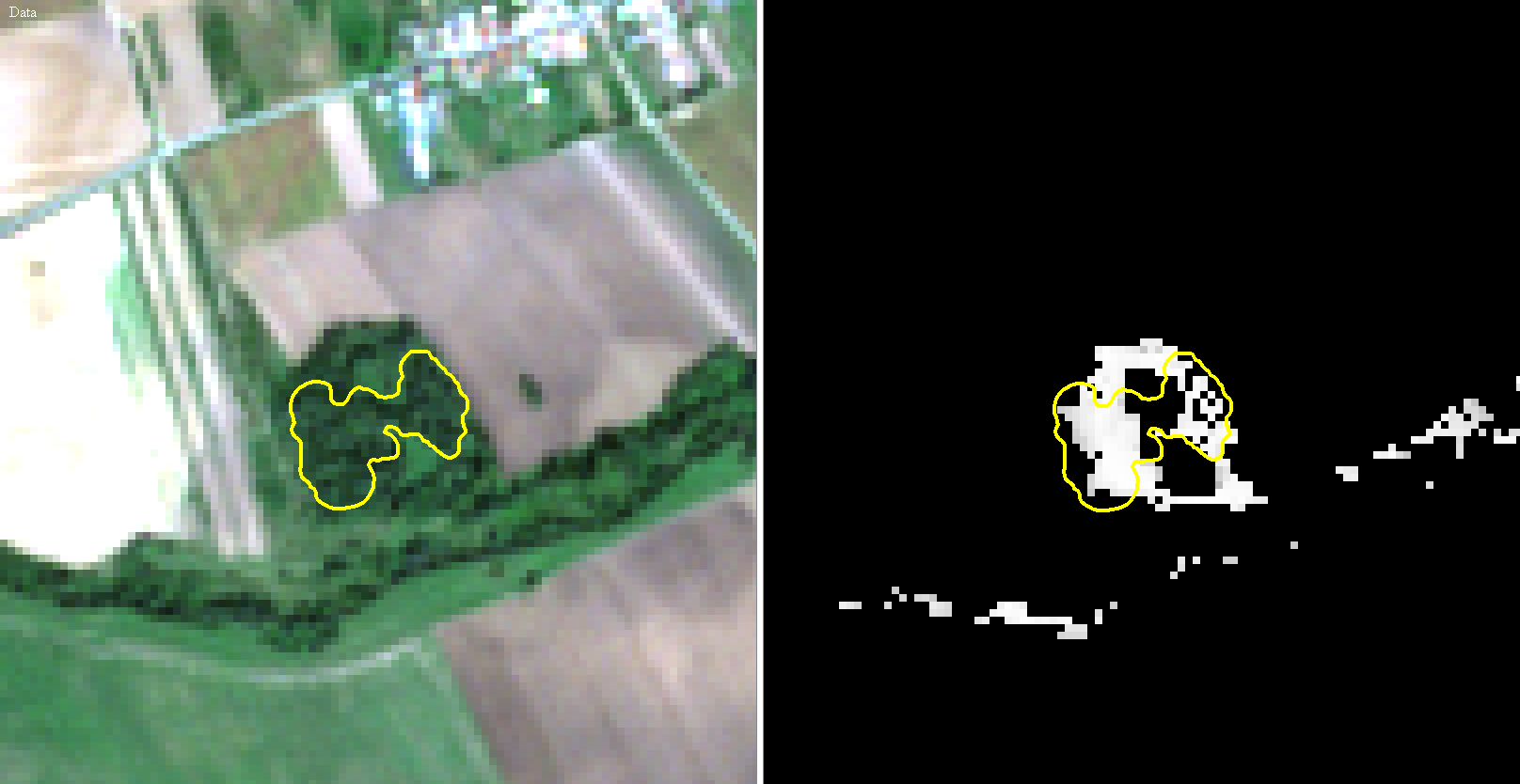} \\
	~\includegraphics[width = 0.48\linewidth]{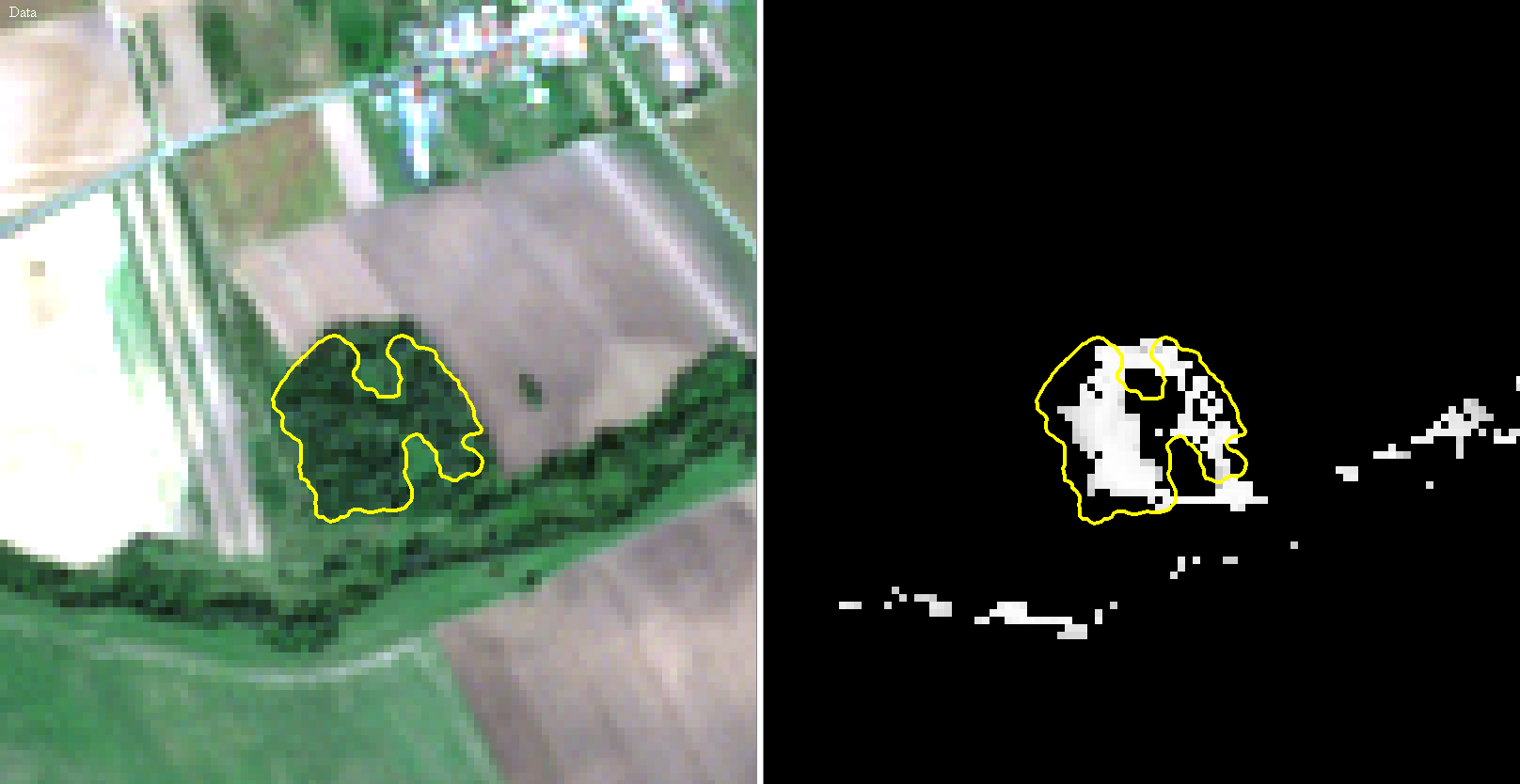}~~
	\includegraphics[width = 0.48\linewidth]{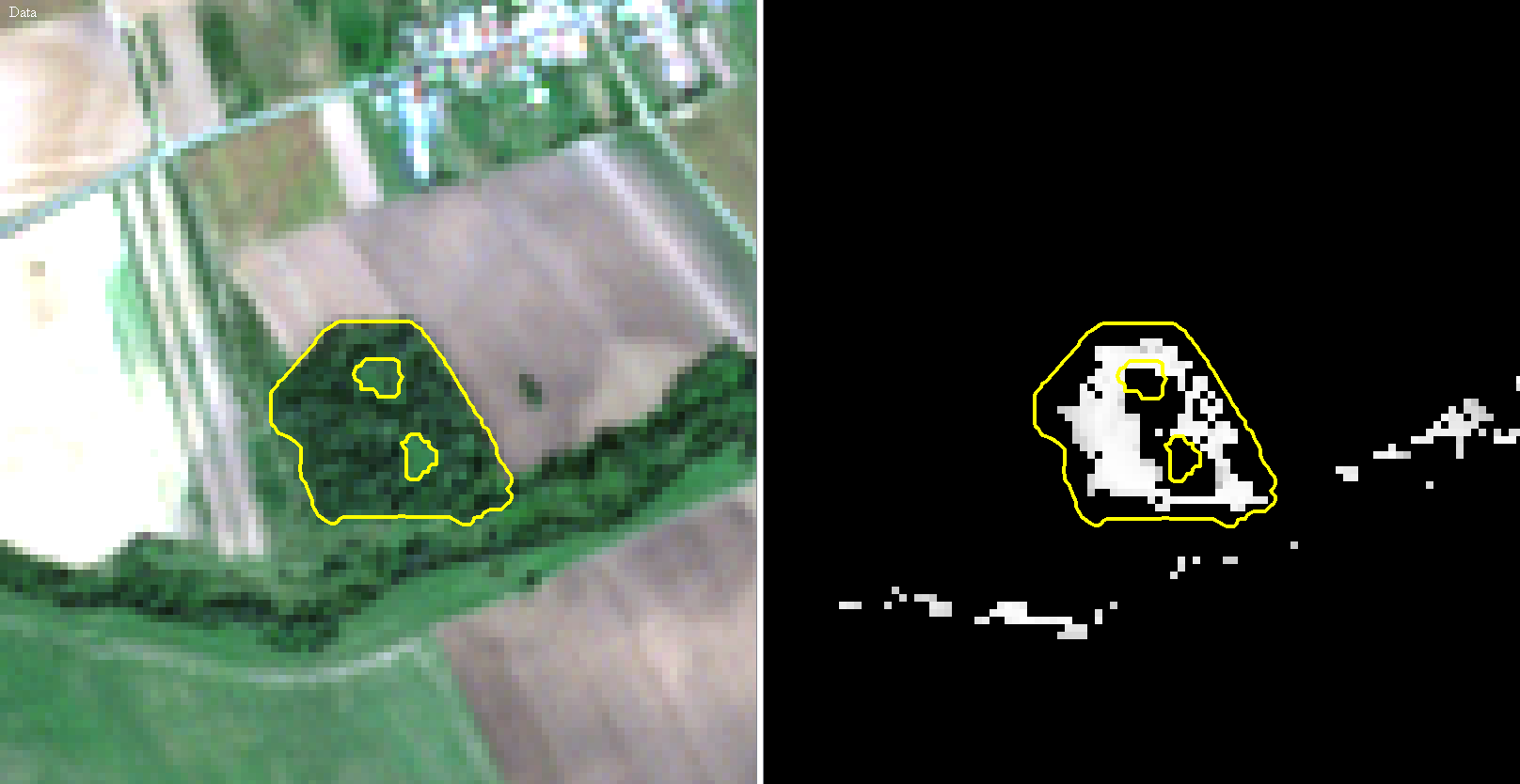}
	\caption{The area around the Rimava river in south of Central Slovakia, the Sentinel-2 image (left) and relevancy map (right) with subsequent automatic segmentation of newly found {\rm 91E0} habitat area.}
	\label{fig:RS}
 \end{center}
\end{figure}

\begin{figure}
 \begin{center}
	\includegraphics[width = 0.48\linewidth]{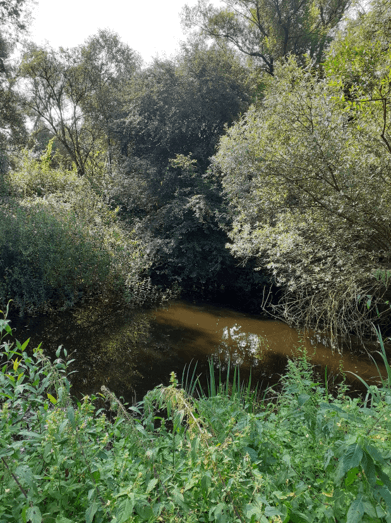}~~
	\includegraphics[width = 0.48\linewidth]{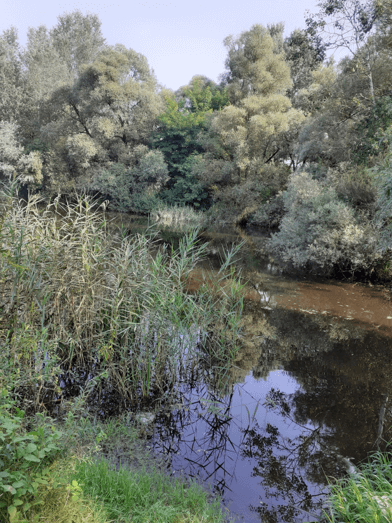}
	\caption{Newly discovered area of {\rm 91E0} protected habitat around the Rimava river, photographs.}
	\label{fig:RS_real}
 \end{center}
\end{figure}

During the exploration of the relevancy maps created for the southern area of the Central Slovakia region, we found that the natural network assigns high 91E0 relevancy to the area around the Rimava River close to the village of Dubovec (Fig. \ref{fig:RS} right). The discovered area has never been identified as a target habitat, and no databases contain such information (the Slovak vegetation database or the database of the State Nature Conservancy of the Slovak Republic where all currently known areas of Natura 2000 habitats are collected). We segmented the area by applying the automatic segmentation method \cite{segAuto}, as shown at the bottom right in Fig. \ref{fig:RS}, to obtain the final segmentation result. The computed mean relevancy inside the segmented area is equal to $0.6079$, and it again indicates a possible new appearance of the 91E0 habitat. Finally, the botanists went into the field and confirmed that the newly discovered area is classified as 91E0 Natura 2000 habitat, as shown in photos of the area in Fig. \ref{fig:RS_real}.
	
\subsection{Exploration of the alluvial forests along the Danube River in the Central and South Europe}
During the field exploration of the Danube River floodplains, the appearance and positions of the habitats were sampled using vegetation plots. We use the natural network to explore the alluvial forests along the Danube River and we made a qualitative and quantitative comparison of the vegetation plots and the relevancy map.

\begin{figure}
 \begin{center}
	\includegraphics[width = 0.48\linewidth]{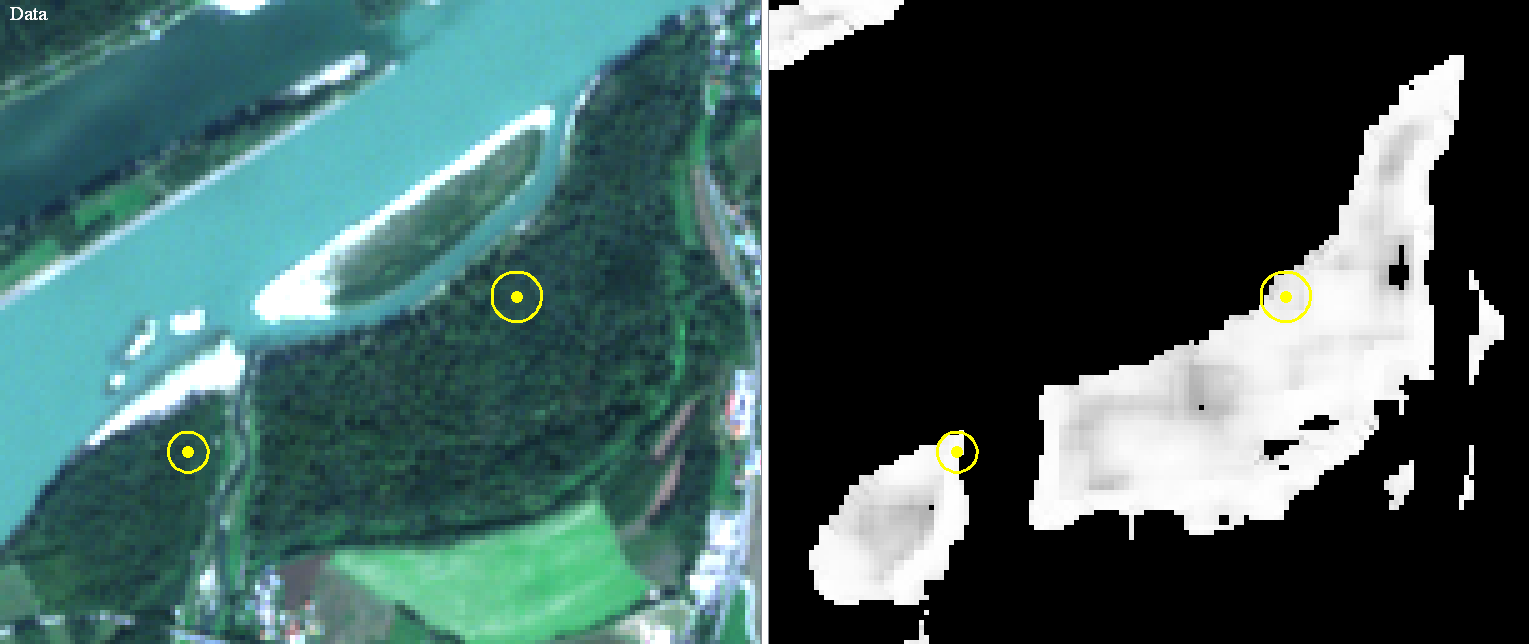}~~
	\includegraphics[width = 0.48\linewidth]{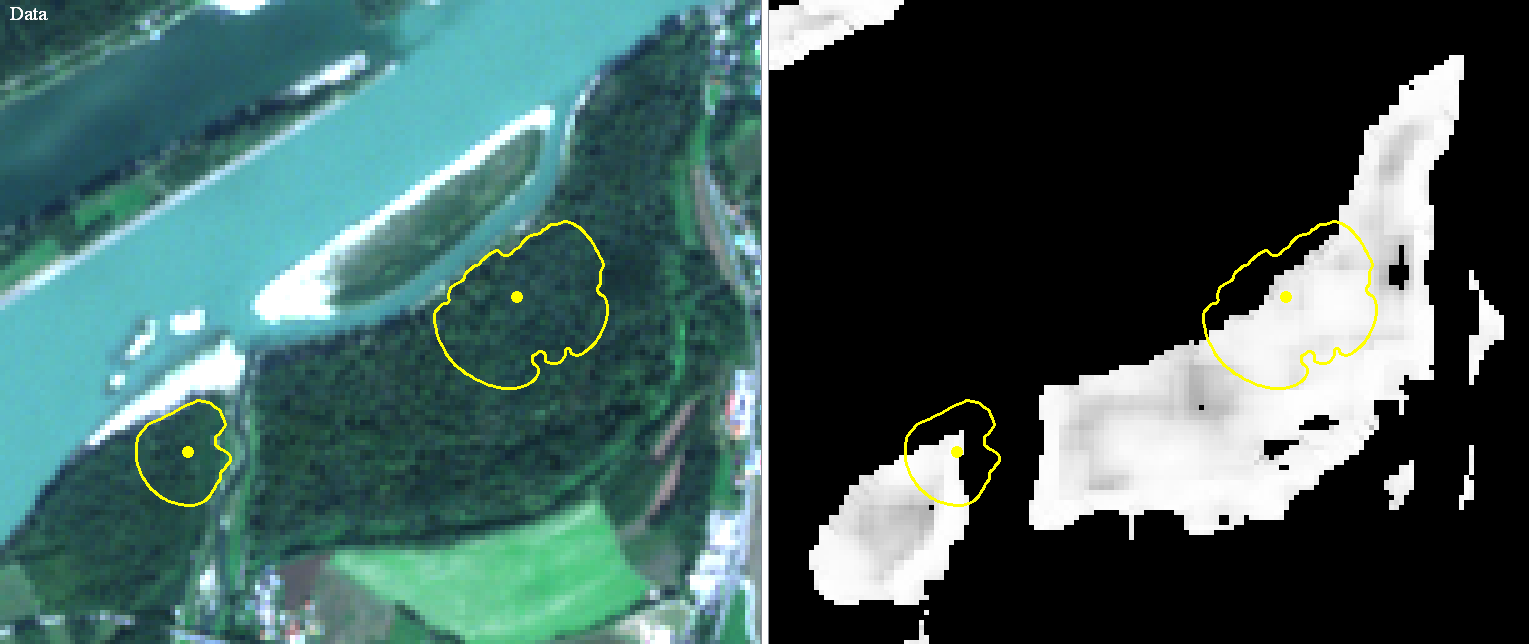} \\
	~\includegraphics[width = 0.48\linewidth]{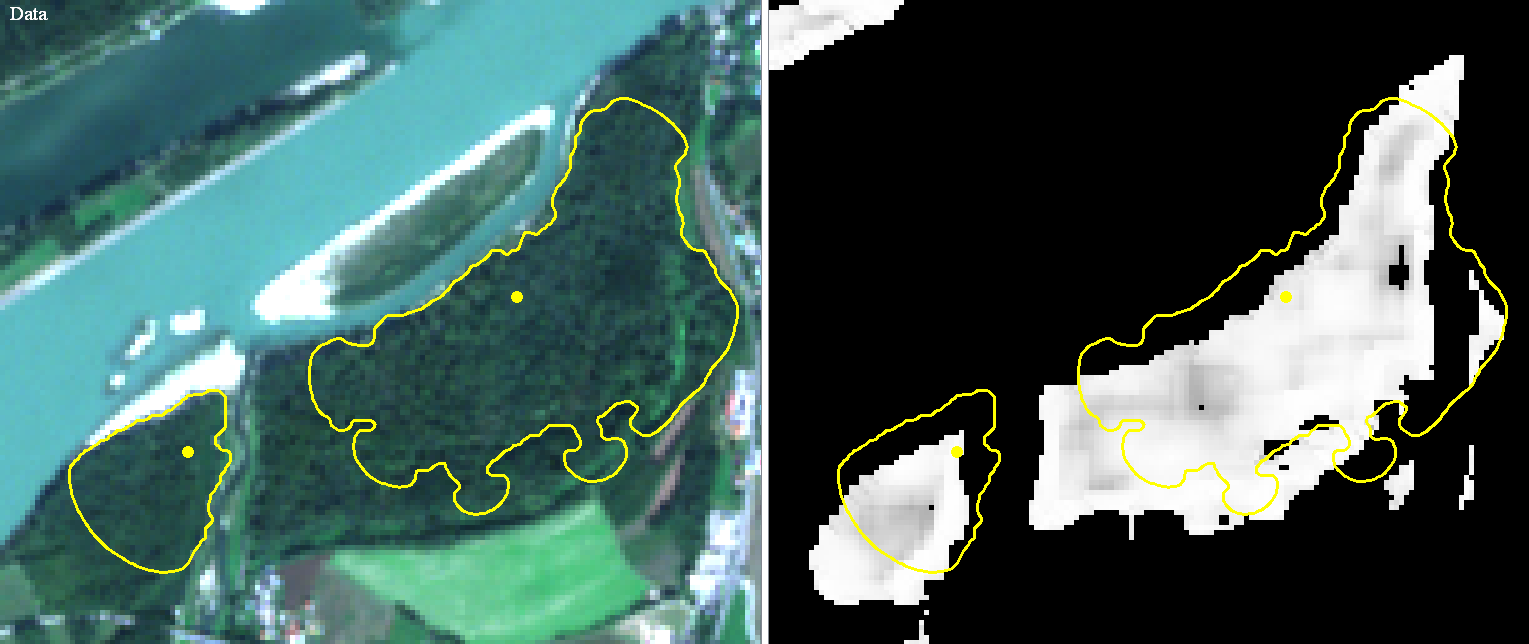}~~
	\includegraphics[width = 0.48\linewidth]{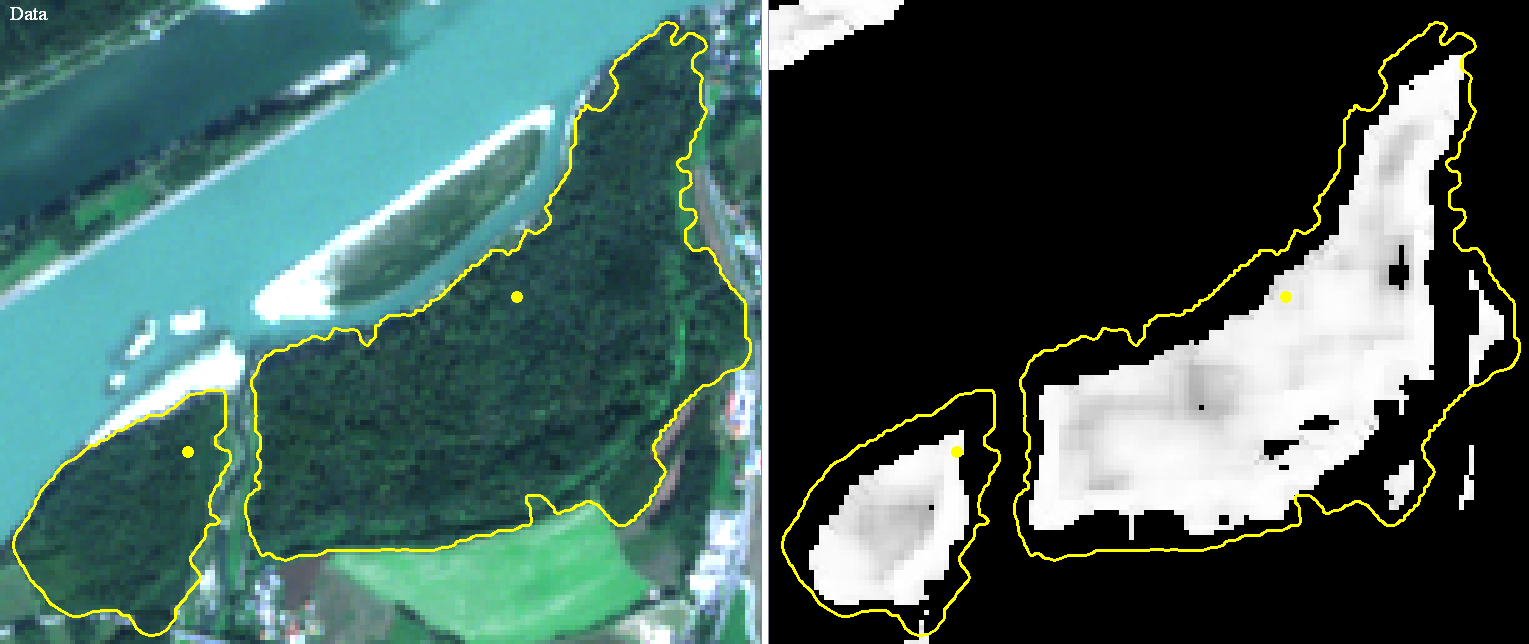}
	\caption{The area on the Danube river bank in Austria around Linz, the Sentinel-2 image (left) and relevancy map (right) with subsequent automatic segmentation of two {\rm 91E0} habitat area.}
	\label{fig:linz}
 \end{center}
\end{figure}

The floodplain forests occur on the banks of the Danube River from Central to South Europe. A large area of these forests is situated in Upper Austria around Linz city, where softwood floodplain forests (91E0 habitat) form a mosaic with monodominant plantations of Canadian poplars or maples. Vegetation plots were sampled in this area, and they were used to verify the relevancy map. All of the vegetation plots fit inside the regions identified by the relevancy map.  Fig. \ref{fig:linz} depicts two of them, situated inside the areas, which the natural network denotes as 91E0 habitat areas. Automatic segmentation was applied to these areas, and the result is also depicted in Fig. \ref{fig:linz}. The final curves surround the areas with high 91E0 relevancy. The 91E0 mean relevancy of the areas was equal to $ 0.7226 $ for the smaller area and $ 0.7318 $ for the larger area.

\begin{figure}
 \begin{center}
	\includegraphics[width = 0.48\linewidth]{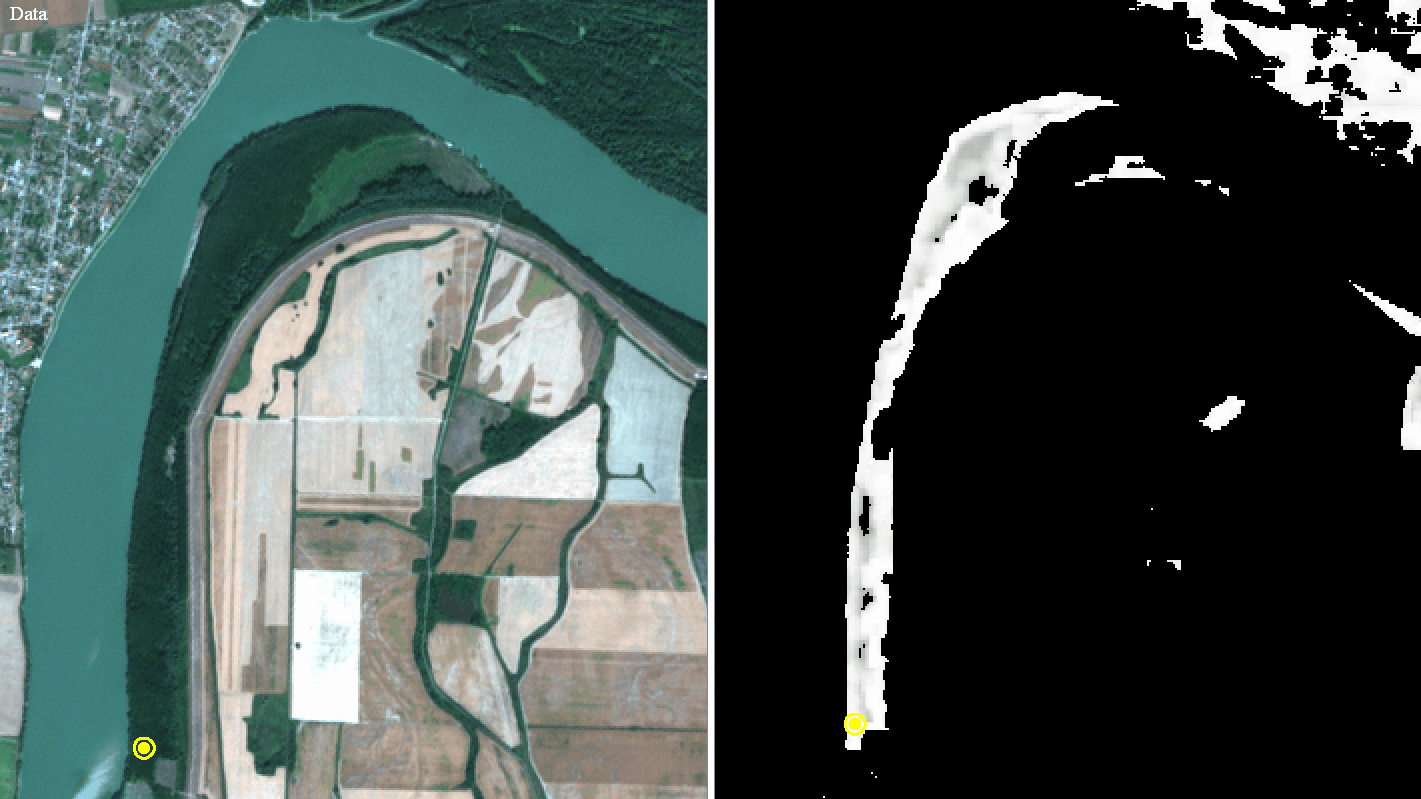}~~
	\includegraphics[width = 0.48\linewidth]{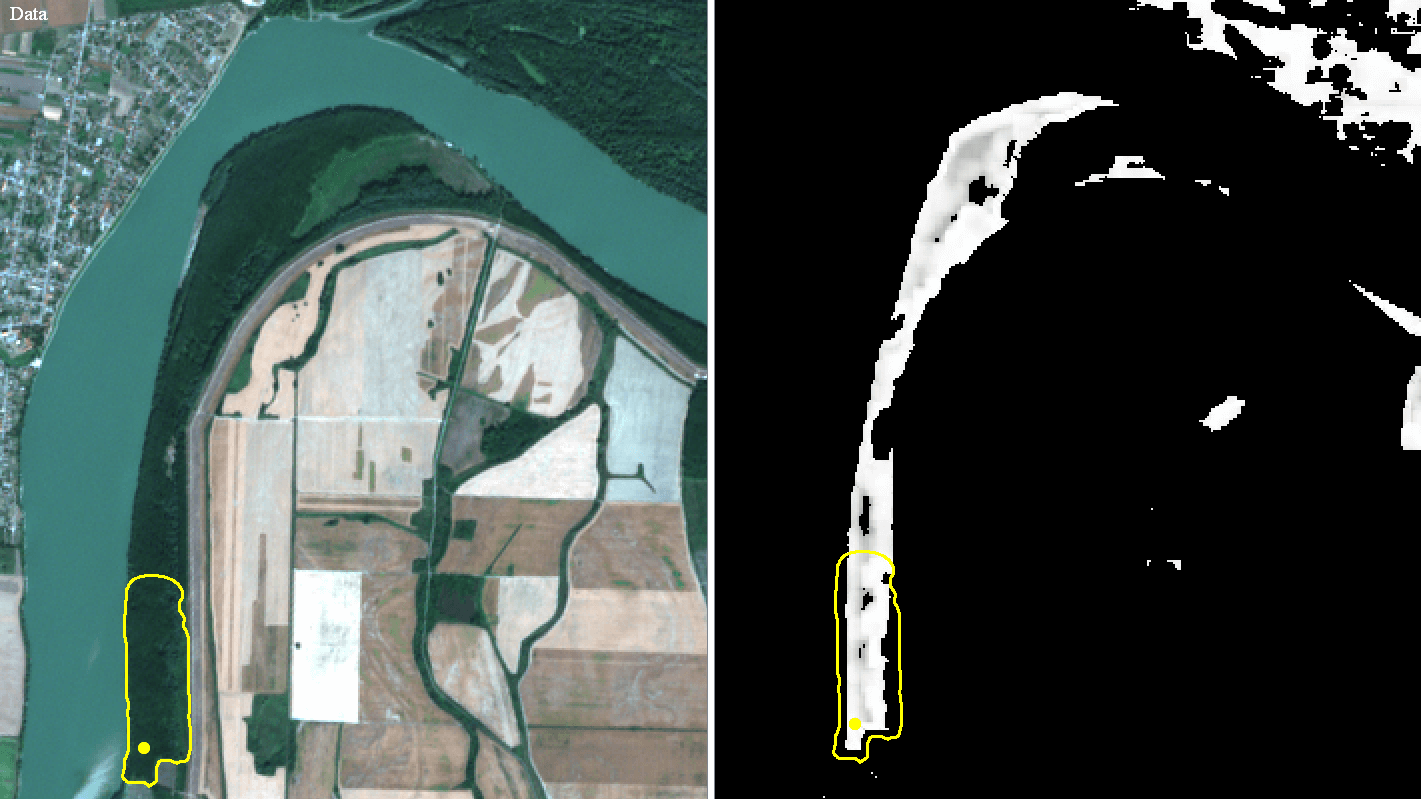} \\
	~\includegraphics[width = 0.48\linewidth]{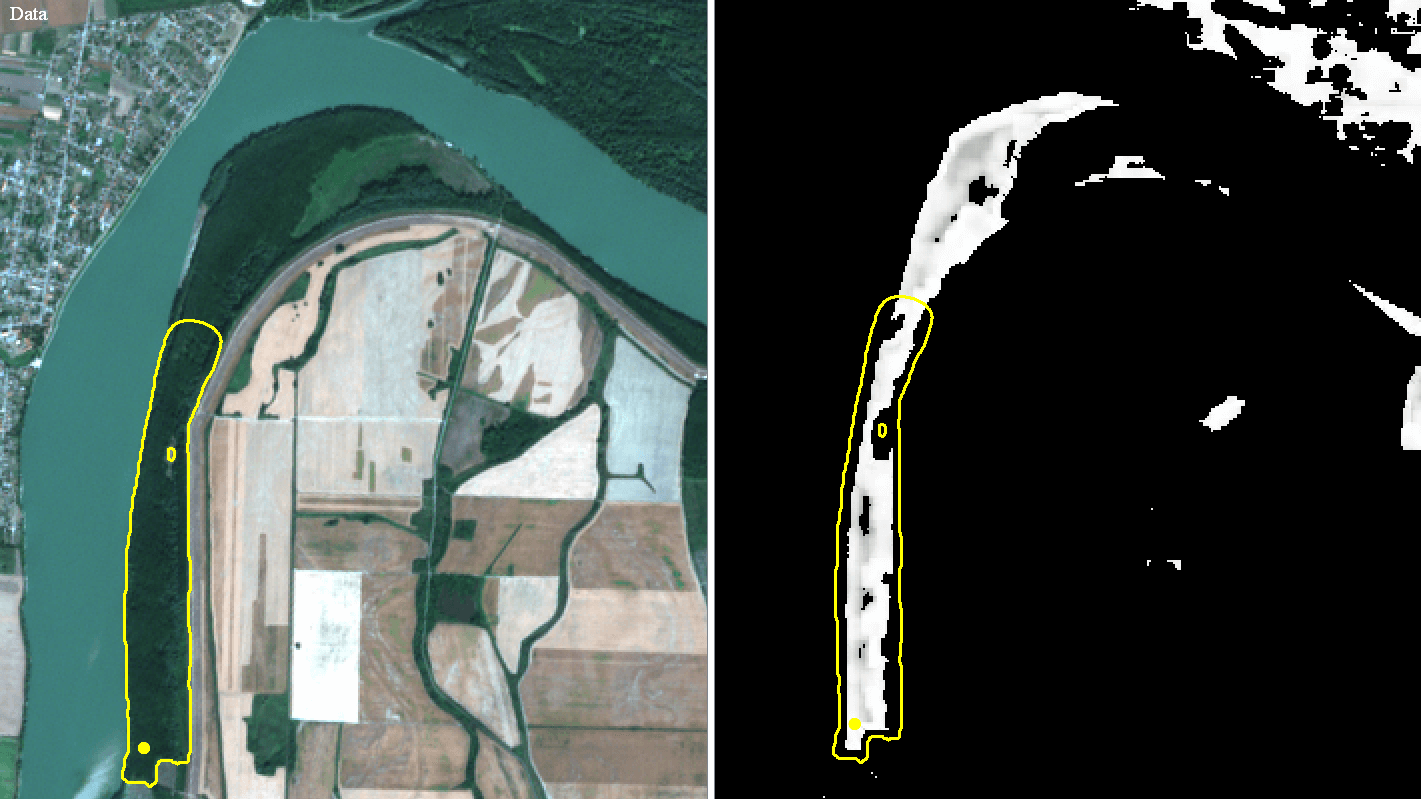}~~
	\includegraphics[width = 0.48\linewidth]{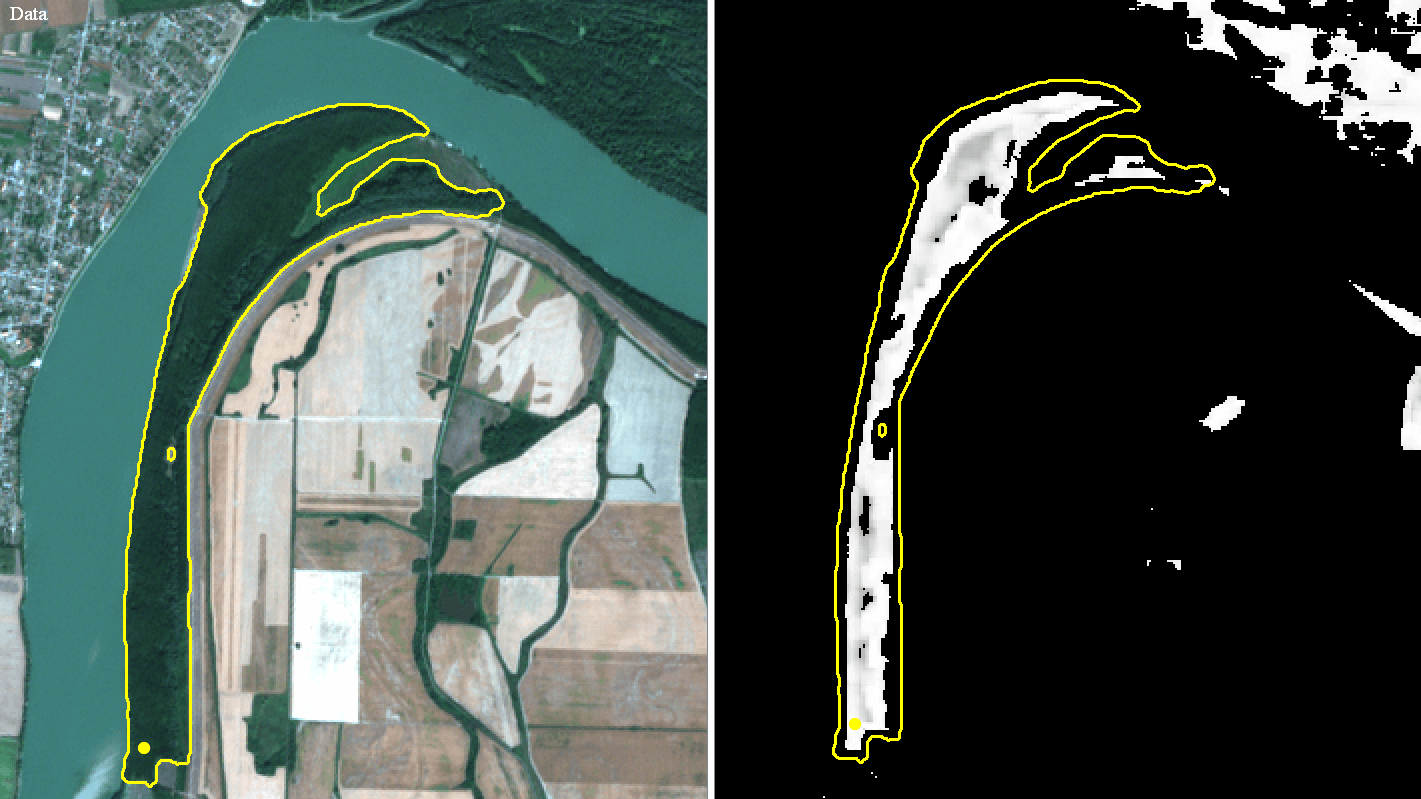}
	\caption{The area around the Danube river in Western Vojvodina, the Sentinel-2 image (left) and relevancy map (right) with subsequent automatic segmentation of {\rm 91E0} habitat area.}
	\label{fig:vajska}
 \end{center}
\end{figure}

\begin{figure}
 \begin{center}
	\includegraphics[width = 0.9\linewidth]{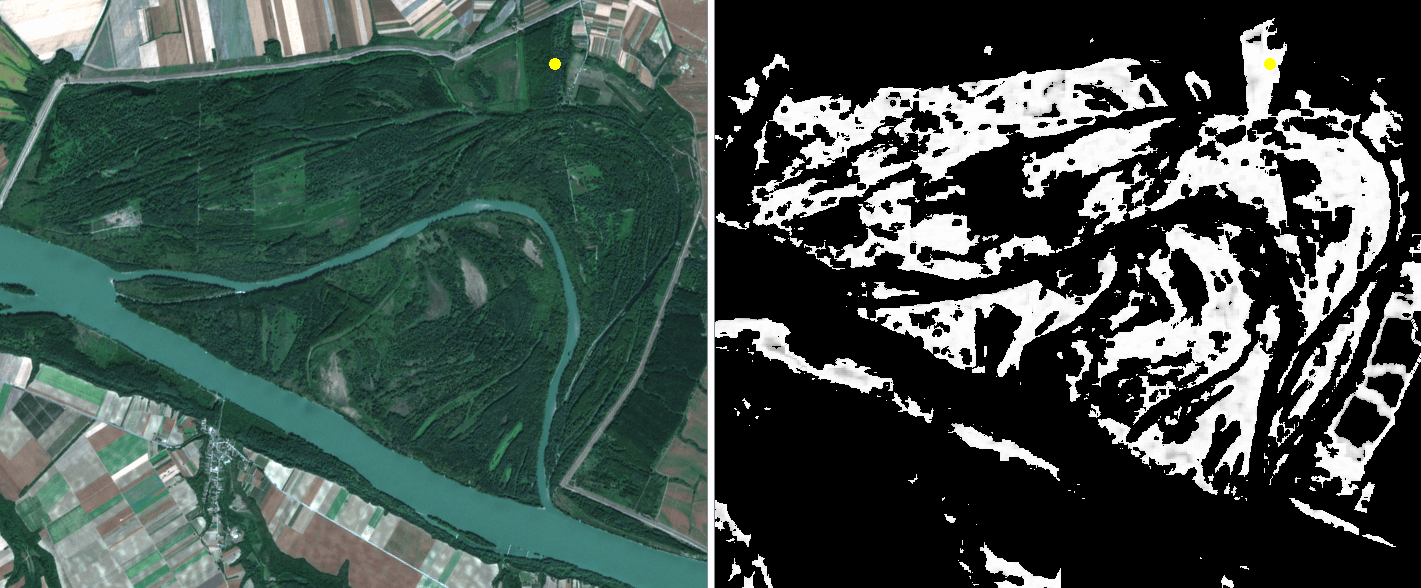} \\
	\includegraphics[width = 0.9\linewidth]{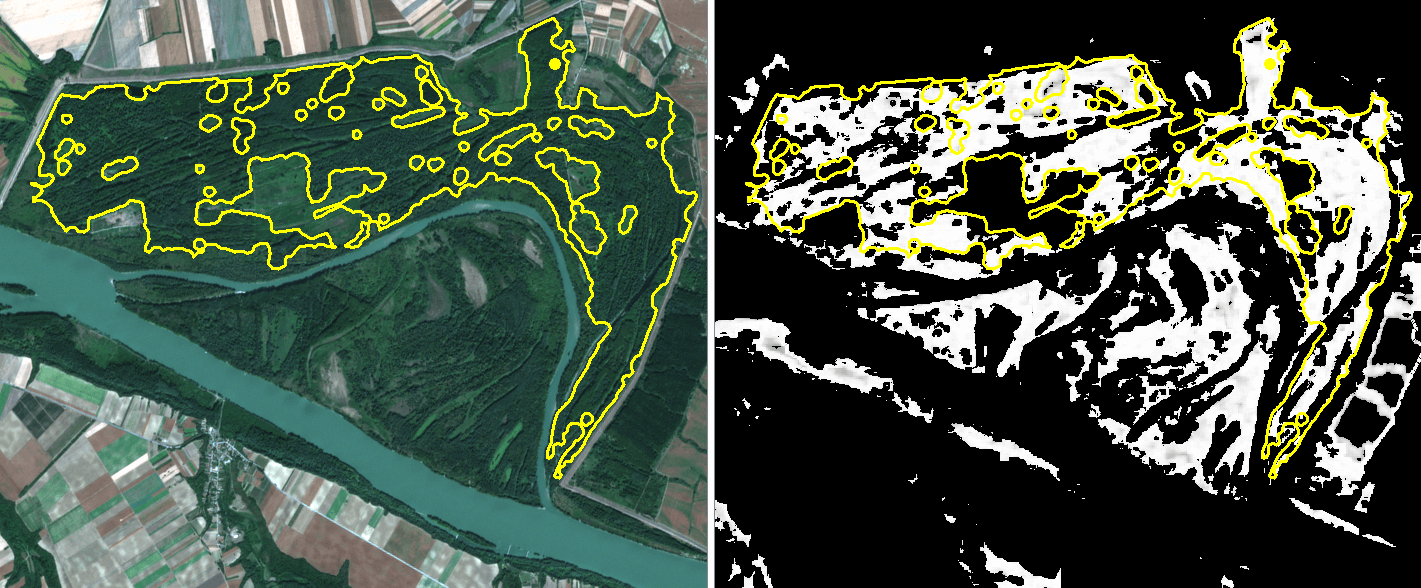}
	\caption{The special nature reservation "Karadordevo" around the Danube river in Western Vojvodina, the Sentinel-2 image (left) and relevancy map (right) with automatic segmented area of {\rm 91E0} habitat.}
	\label{fig:kara}
 \end{center}
\end{figure}

The relevancy map was created, and vegetation plots were also sampled in the Danube River alluvia in Vojvodina (North Serbia). The identification of 91E0 habitat areas was slightly less successful in this area. Two plots were not recognized due to the presence of the moistest forest types dominated by sparse willows and absenting poplars that are not so widespread in the upper parts of the Danube alluvia where the natural network was trained. One habitat area was not found because the plot occurs along a very thin forest line, which is explainable by the Sentinel-2’s resolution and the mechanism of creating the relevancy map. All other plots were identified correctly, and two examples of habitat areas in the relevancy map are shown in Fig. \ref{fig:vajska} and Fig. \ref{fig:kara}, respectively. The first example is situated on the left bank of the Danube River near the village of Vajska. This vegetation plot was used as a starting point for the automatic segmentation algorithm, as shown in Fig. \ref{fig:vajska}. This thin habitat area was segmented with a 91E0 mean relevancy equal to $0.7091$. 

A complex habitat system from the Karadordevo nature reservation, identified in the relevancy map and segmented by automatic segmentation, is presented in Fig. \ref{fig:kara}. The starting point for the segmentation was the vegetation plot located on the north border of this area. The entire area of the 91E0 habitat was successfully segmented while the mosaic of wet meadows and young clear-cuts inside the habitat area was excluded. The mean relevancy of the segmented area is $0.6032$, which means correct classification within the 91E0 habitat.

\begin{figure}
 \begin{center}
	\includegraphics[width = 0.9\linewidth]{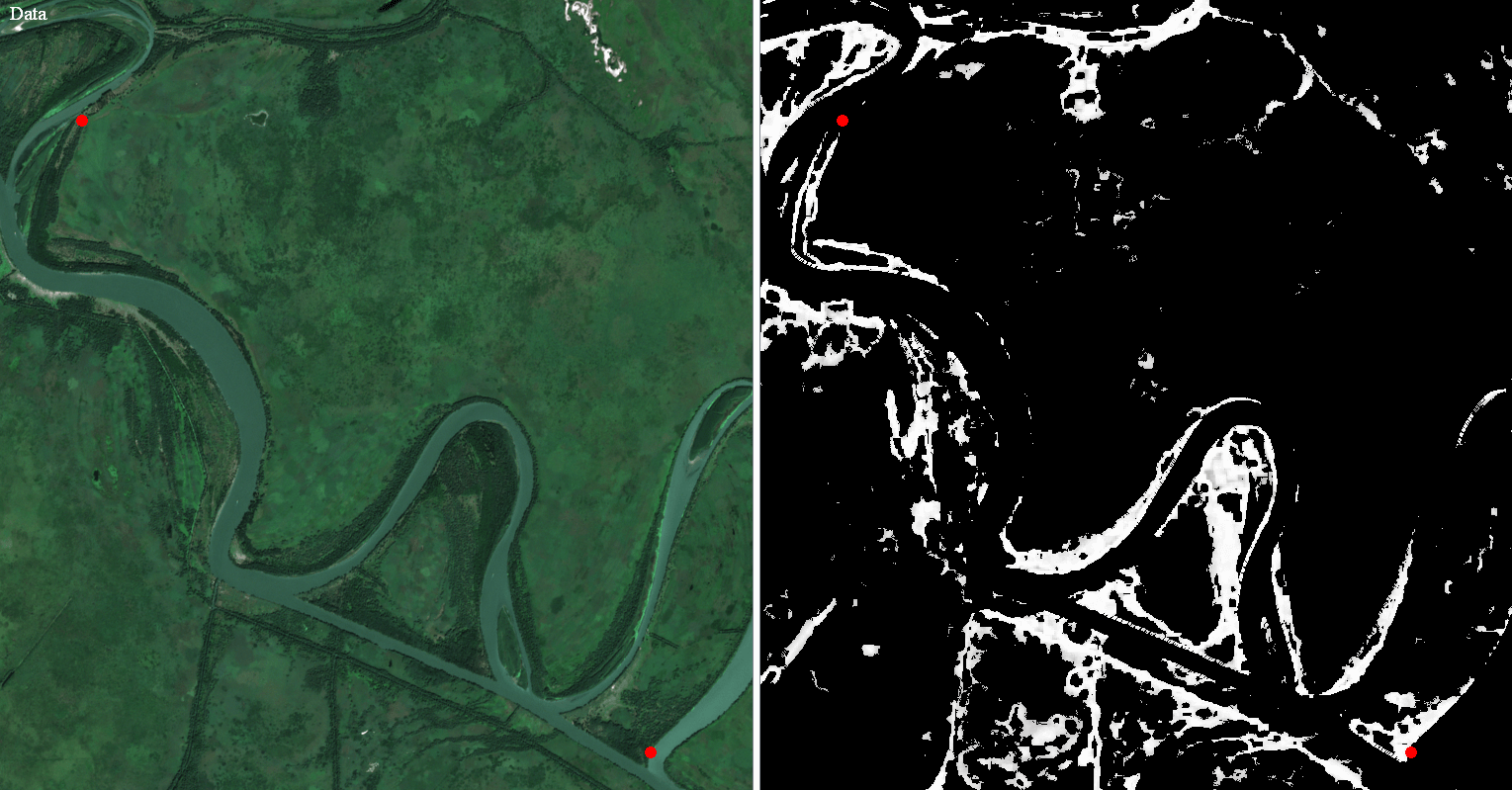} \\
	~\includegraphics[width = 0.9\linewidth]{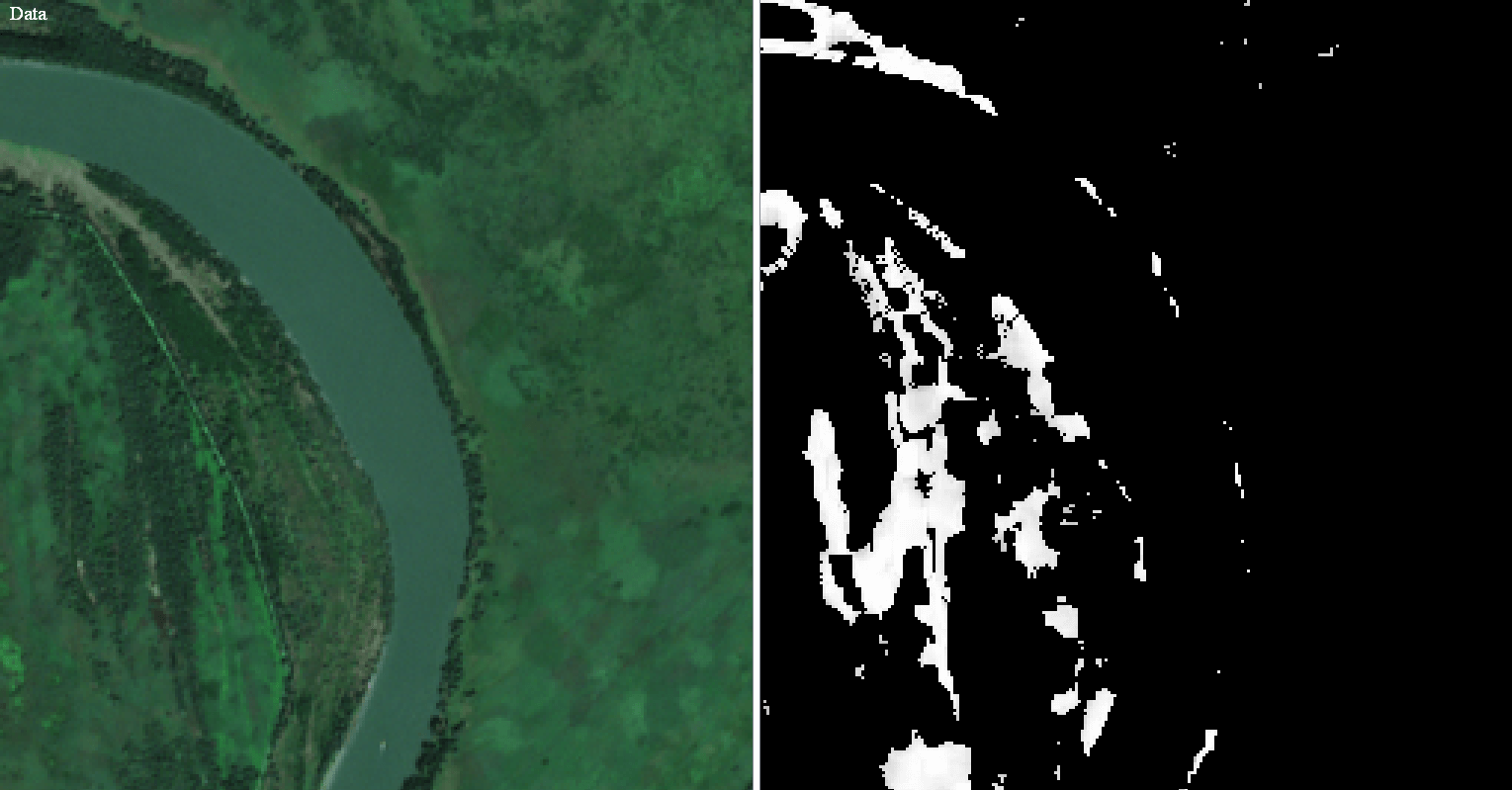}
	\caption{The delta of the Danube river in Romania, the Sentinel-2 images (left) and details of the relevancy map (right). By the red color points, the 91E0 habitat localities, with the photographs given in Fig. \ref{fig:DJ}, are indicated.}
	\label{fig:deltaJuh}
 \end{center}
\end{figure}

Another studied area was the delta of the Danube River in Romania. It is difficult to explore the region due to the complicated accessibility of most of the area. The character of the landscape means that localities can only be reached by small boats; thus, the importance of creating relevancy maps of such areas is significant. Fig. \ref{fig:deltaJuh} plots the Saint George Branch of the Danube River, and the relevancy map corresponds very well with the vegetation sampled by plots. In the detail of the Danube delta's branch relevancy map, we can clearly see bright-colored 91E0 habitat and black-colored surrounding swamps and other ecosystems, which indicates correct classification.


\begin{figure}
 \begin{center}
	\includegraphics[width = 0.9\linewidth]{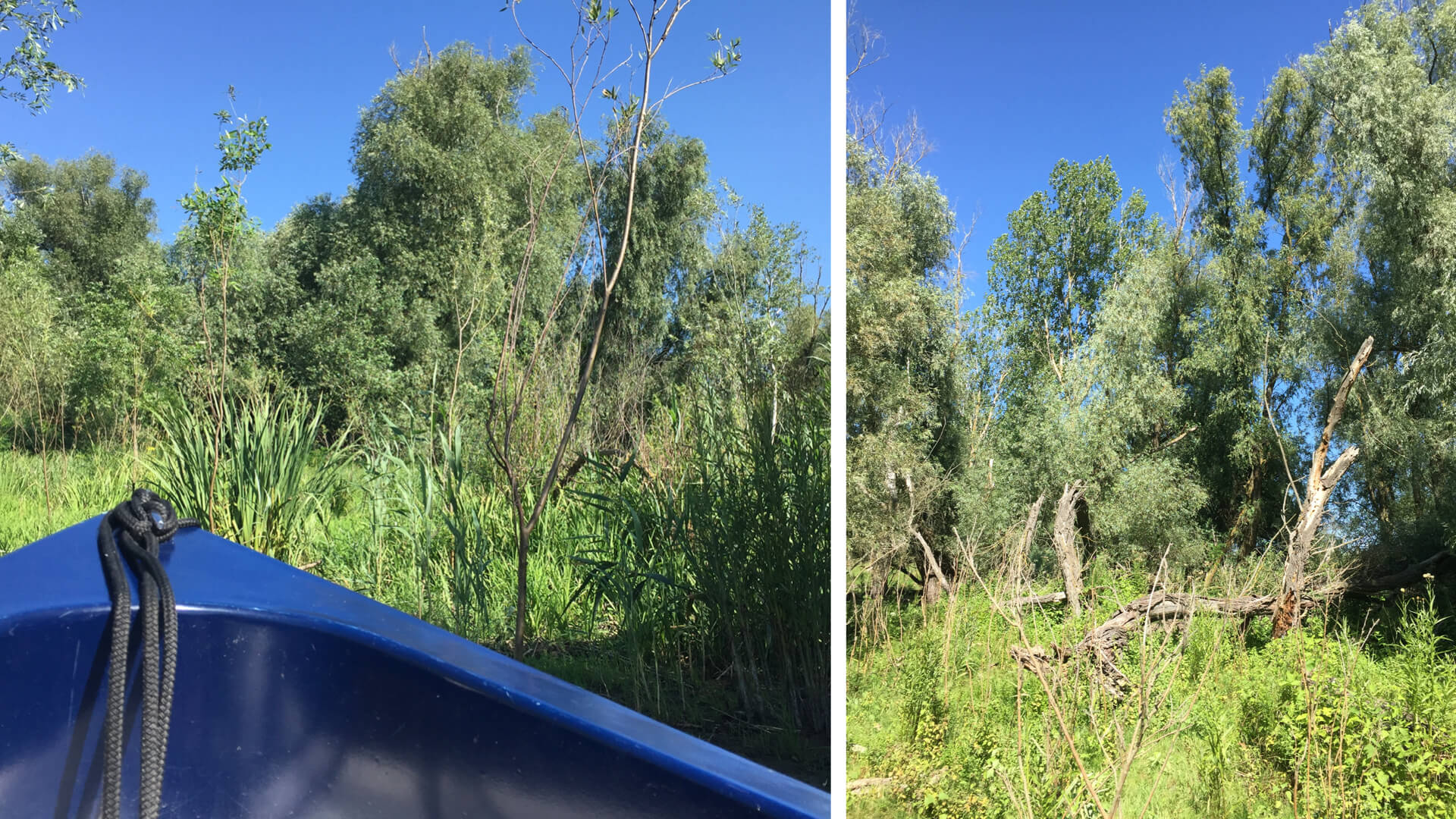}
	\caption{The delta of the Danube river in Romania, photographs.}
	\label{fig:DJ}
 \end{center}
\end{figure}

\section{Conclusions} The presented results give new perspectives for diverse science disciplines - mathematical modelling, vegetation science and ecology, remote sensing, nature conservation and mapping, and subsequent delivery of ecosystem services. The identification of plant communities on the scale of, e.~g., Natura 2000 habitats, using remote sensing has been an open challenge for field ecologists over the last decades. Rapidly developing remote sensing techniques and different data-mining approaches were implemented to monitor land surface types such as agricultural land, water bodies, abandonment land, natural and plantation forests of different types, meadows with various management practices, built-up areas, etc, see e.g. \cite{Franklin_etal.2001, Fagan_etal.2015, Laurin_etal.2016, Noi, Erinjery_etal.2018, Cheng&Wang2019, Wasnievski_etal.2020}. However, due to the complicated character of target nature phenomena, it was not possible to reach the detailed scale of Natura 2000 habitats by using existing methodologies and satellite data. Thus, we developed the novel method - the natural numerical network - for that purpose. We introduced its definition in the form of the discretized forward-backward diffusion equation and applied it to the classification of Natura 2000 forest habitats. We also introduced the new concept of relevancy maps, which allow the automatic recognition of Natura 2000 habitat areas in remote sensing data provided by Sentinel-2 optical information and finding new appearances of protected habitats via the satellite data. We are not aware of any other method that would identify and explore the Natura 2000 habitats or similarly detailed plant communities on such an exact level and accuracy by using the remote sensing data.

\section{Acknowledgement}
This work was supported by the grants APVV-16-0431, APVV-19-0460 and ESA contracts 4000122575/17/NL/SC and 4000133101/20/NL/SC.

\bibliography{bibliography}

\end{document}